\documentclass[moor,nonblindrev]{informs1} 

\usepackage{natbib}
 \NatBibNumeric
 \def\bibfont{\small}%
 \bibpunct[, ]{[}{]}{,}{n}{}{,}%

\usepackage[colorlinks=true,breaklinks=true,bookmarks=true,urlcolor=blue,
     citecolor=blue,linkcolor=blue,bookmarksopen=false,draft=false]{hyperref}

\def\EMAIL#1{\href{mailto:#1}{#1}}

\usepackage[inline]{enumitem}
\usepackage{tikz}
\usepackage{subcaption}			
\usepackage{multirow}			
\usepackage[ruled]{algorithm2e}
\usepackage{multicol} 	

\newcommand{\B}{\{0,1\}}
\newcommand{\Z}{\mathbb{Z}}
\newcommand{\RR}{\mathcal{R}}
\newcommand{\R}{\mathbb{R}}
\newcommand{\mvec}[1]{\boldsymbol{#1}}
\DeclareMathOperator\supp{supp}

\TheoremsNumberedThrough     

\EquationsNumberedThrough    

\MANUSCRIPTNO{} 

\begin{document}


\RUNAUTHOR{G\"unl\"uk, Hauser and Kov\'acs}

\RUNTITLE{Binary Matrix Factorisation and Completion via Integer Programming}

\TITLE{Binary Matrix Factorisation and Completion via Integer Programming}

\ARTICLEAUTHORS{%
\AUTHOR{Oktay G\"unl\"uk}
\AFF{Cornell University, \EMAIL{ong5@cornell.edu}}
\AUTHOR{Raphael A. Hauser, R\'eka \'A. Kov\'acs}
\AFF{University of Oxford, The Alan Turing Institute, \EMAIL{hauser@maths.ox.ac.uk}, \EMAIL{reka.kovacs@maths.ox.ac.uk}}
} 

\ABSTRACT{%
	Binary matrix factorisation is an essential tool for identifying discrete patterns in binary data.
In this paper we consider the rank-$k$ binary matrix factorisation problem ($k$-BMF) under Boolean arithmetic: we are given an $n\times m$ binary matrix $X$ with possibly missing entries and need to find two binary matrices $A$ and $B$ of dimension $n\times k$ and $k\times m$ respectively, which minimise the distance between $X$ and the Boolean product of $A$ and $B$ in the squared Frobenius distance.
We present a compact and two  exponential size integer programs (IPs) for $k$-BMF and show that the compact IP has a weak LP relaxation, while the exponential size IPs have a stronger equivalent LP relaxation. We introduce a new objective function, which differs from the traditional squared Frobenius objective in attributing a weight to zero entries of the input matrix that is proportional to the number of times the zero is erroneously covered in a rank-$k$ factorisation. 
For one of the exponential size IPs we describe a computational approach based on column generation.
Experimental results on synthetic and real word datasets suggest that our integer programming approach is competitive against available methods for $k$-BMF and provides accurate low-error factorisations.
}%


\MSCCLASS{90C10}
\ORMSCLASS{Integer Programming}
\HISTORY{}

\maketitle

%


\section{Introduction.}

For a given binary matrix $X\in\B^{n\times m}$ and a fixed positive integer $k$, the rank-$k$ binary matrix factorisation problem ($k$-BMF) is concerned with finding two matrices $A\in\B^{n\times k}$, $B\in \B^{k\times m}$ such that the product of $A$ and $B$ is a binary matrix closest to $X$ in the squared Frobenius norm.
One can define different variants of this problem depending on the underlying arithmetic used when computing the product of the matrices.
In this paper we focus on solving $k$-BMF under Boolean arithmetic where the product of the binary matrices $A$ and $B$ is computed by 
$(i)$ interpreting $0$s as false and $1$s as true,
and 
$(ii)$ using logical disjunction ($\vee$) in place of addition and logical conjunction ($\wedge$) in place of multiplication.
Observe that Boolean multiplication ($\wedge$) coincides with standard multiplication on binary input, hence we adopt the  notation $a\,b$ in place of $a\wedge b$ in the rest of the paper.
We therefore compute the Boolean matrix product of  $A$ and $B$ as:
\begin{equation*}
	Z= A\circ B \iff z_{i j} = \bigvee_{\ell}(a_{i\ell}\, b_{\ell j}). 
\end{equation*}
Note that Boolean matrix multiplication can be equivalently written as  $z_{ij} = \min\{1,\sum_{\ell} a_{i\ell} b_{\ell j}\}$ using standard arithmetic summation.
The problem then becomes computing matrices $A$ and $B$ whose Boolean product $Z$ best approximates the input matrix $X$.

Our motivation for this study comes from data science applications where rows of the matrix $X$ correspond to data points and columns correspond to features.
In these applications low-rank matrix approximation is an essential tool for dimensionality reduction which helps understand the data better by exposing hidden features. 
Many practical datasets contain categorical features which can be represented by a binary data matrix using unary encoding. 
For example, consider a data matrix $X$ below (inspired by \cite{Miettinen:2008:DBP:1442800.1442809}),  
where rows correspond to patients and columns to symptoms, $x_{ij}=1$ indicating patient $i$ presents symptom $j$:
\begin{align}
	\label{example_intoduction}
	X=\begin{bmatrix}		1 & 1 & 0\\		1 & 1 & 1\\		0 & 1 & 1	\end{bmatrix}
	&&
	X = 	A\circ B = 	\begin{bmatrix}		1 & 0 \\		1 & 1 \\		0 & 1	\end{bmatrix}
	\circ 	\begin{bmatrix}		1 & 1 & 0 \\		0 & 1 & 1	\end{bmatrix}.
\end{align}
In this example  matrix $A\circ B$ describes $X$ exactly using  $2$ derived features where the rows of $B$ specify how the original features relate to the $2$  derived features, and  the rows of $A$ give the derived features of each patient.
In other words, factor matrix $B$ reveals that there are 2 underlying diseases that cause the observed symptoms: Disease $\alpha$ is causing symptoms 1 and 2, and disease $\beta$ is causing symptoms 2 and 3. Matrix $A$ reveals that patient 1 has disease $\alpha$, patient 3 has $\beta$ and patient 2 has both.

We note that it is also possible to use classical methods such as singular value decomposition (SVD) \cite{SVD} or  non-negative matrix factorisation (NMF) \cite{NMF}  to obtain low-rank approximations of $X$ but the resulting factor matrices or their product would typically not be binary unlike  BMF \cite{Miettinen:2006:PKDD}.
To demonstrate this we next give the best rank-$2$  SVD and NMF approximations of the matrix $X$ in \eqref{example_intoduction}, respectively:
\begin{align}
	\label{example_intoduction_part2}
	X 	\approx
	\begin{bmatrix}		1.21&0.71\\		1.21 &0.00\\		1.21 &-0.71\\	\end{bmatrix}
	\begin{bmatrix}		0.00&0.71&0.50\\		0.71& 0.00&-0.71	\end{bmatrix},
	&&
	X 	\approx
	\begin{bmatrix}		1.36 &0.09\\		1.05 &1.02\\		0.13 &1.34	\end{bmatrix}
	\begin{bmatrix}		0.80 &0.58 &0.01\\		0.00     &     0.57 &0.81	\end{bmatrix}.
\end{align}
Note that neither of these rank-2  approximations provide a clear interpretation. 
The rank-2 NMF of $X$ suggests that symptom 2 presents with lower intensity in both $\alpha$ and $\beta$, an erroneous conclusion (caused by patient 2) that could not have been learned from data $X$ which is of ``on/off'' type.

We note that in addition to healthcare applications, BMF-derived features of data have also been shown to be interpretable in biclustering gene expression datasets \cite{Zhang:2007}, role based access control \cite{Lu:2008:OBM:1546682.1547186,Lu:2014:9400730720140101} and market basket data clustering \cite{Li:2005}.

\subsection{Complexity and related work.}

The Boolean rank \cite{Monson:1995, Kim:1982} of a binary matrix $X$ is defined to be the smallest integer $r$ for which there exist binary matrices $A$ and $B$  such that $X=A\circ B$.
In an equivalent definition, the Boolean rank of $X$ is the minimum value of $r$ for which it is possible to factor $X$ into a Boolean combination of $r$ rank-$1$ binary matrices $$X = \bigvee_{\ell=1}^r \mvec{a}_\ell \, \mvec{b}^\top_\ell$$ for 
$\mvec{a}_\ell \in\B^{n},\mvec{b}_\ell \in \B^m$. Occasionally, the Boolean rank is also referred to as the rectangle cover number, and rank-$1$ binary matrices $\mvec{a}_\ell \mvec{b}^\top_\ell$ are called rectangle matrices or simply rectangles \cite{Conforti:2014:IP:2765770}.

Interpreting $X$ as the node-node incidence matrix of a bipartite graph $G(X)$ with $n$ vertices on the left and $m$ vertices on the right, the problem of computing the Boolean rank of $X$ is in one-to-one correspondence with finding a minimum edge covering of $G(X)$ by complete bipartite subgraphs (bicliques)\cite{Monson:1995}. 
Since the biclique cover problem is NP-hard \cite[Theorem 8.1]{Orlin:1977},\cite[Problem GT18]{Garey:1979:CIG:578533}, 
and hard to approximate \cite{Simon:1990OnAS, Chalermsook:2014}, 
computing the Boolean rank is hard as well. 
Finding an optimal rank-$k$ binary factorisation of $X$ under Boolean arithmetic has a graphic interpretation of 
minimizing the number of errors in an approximate covering of $G(X)$ by $k$ bicliques which are allowed to overlap. 
In the rank-$1$ case the Boolean arithmetic coincides with standard arithmetic and $1$-BMF can be interpreted as computing a maximum weight biclique on the complete bipartite graph $K_{n,m}$ whose edges that are in $G(X)$ have weight $1$ and others weight $-1$. The maximum edge biclique problem with edge weights in $\{-1,1\}$ is NP-hard \cite{Gillis:2018}, 
hence even the computation of a rank-$1$ BMF is computationally challenging.

Due to the hardness results, the majority of methods developed for BMF rely on heuristics. 
The earliest heuristic for BMF, Proximus \cite{Koyuturk:2002:ATA:645806.670310, Koyuturk:2003:PFA:956750.956770}, computes BMF under standard arithmetic using a recursive partitioning idea and computing $1$-BMF at each step. 
Since Proximus, much research has focused on computing efficient and accurate methods for $1$-BMF.
\cite{Shen:2009:MDP:1557019.1557103} proposes an integer program (IP) for $1$-BMF and several relaxations of it, one of which leads to a $2$-approximation,
while \cite{Shi:2014:6899417} provides a rounding based $2$-approximation.
In \cite{Beckerleg:2020} an extension of the Proximus framework is explored which uses the formulations from \cite{Shen:2009:MDP:1557019.1557103} to compute $1$-BMF at each step. 
$k$-BMF under Boolean arithmetic is explicitly introduced in \cite{Miettinen:2006:PKDD, Miettinen:2008:DBP:1442800.1442809}, along with a heuristic called ASSO, which is based on an association rule-mining approach. 
ASSO is further improved in \cite{Barahona:2019} into an alternating iterative heuristics. 
Another approach based on an alternating style heuristic is explored in \cite{Zhang:2007} to solve a non-linear unconstrained formulation of $k$-BMF with penalty terms in the objective for non-binary entries.

In \cite{Lu:2008:OBM:1546682.1547186,Lu:2014:9400730720140101} a series of integer programs for $k$-BMF and exact BMF are introduced. These IPs have exponentially many variables and constraints and require an explicit enumeration of the $2^m$ possible binary row vectors for factor matrix $B$. To tackle the exponential explosion of rows considered, a heuristic row generation using association rule mining and subset enumeration is developed. 
An exact linear IP for $k$-BMF with polynomially many variables and constraints is presented in our previous work \cite{Kovacs:2017}. This model uses McCormick envelopes \cite{McCormick:1976} to linearize the quadratic terms coming from the matrix product. 
We note that both of these integer programs for $k$-BMF, as well as any other element-wise models can be naturally applied in the context of rank-$k$ binary matrix completion by simply setting the objective coefficients corresponding to missing entries to $0$.

\subsection{Our contribution.}

In this paper, we present a comprehensive study on integer programming methods for $k$-BMF. We examine three integer programs in detail: our compact formulation introduced in \cite{Kovacs:2017}, the exponential formulation of \cite{Lu:2008:OBM:1546682.1547186} and a new exponential formulation which we introduced in a preliminary version of this paper in \cite{Kovacs:2020}. 
We prove several results about the strength of LP-relaxations of the three formulations and their relative comparison. In addition, we show that the new exponential formulation overcomes several limitations of earlier approaches. In particular, it does not suffer from permutation symmetry and it does not rely on heuristically guided pattern mining.
Moreover,  it has a stronger LP relaxation than that of \cite{Kovacs:2017}. 
On the other hand, our new formulation has an exponential number of variables which we tackle using a column generation approach that effectively searches over this exponential space without  explicit enumeration, unlike the complete enumeration used for the exponential size model of \cite{Lu:2008:OBM:1546682.1547186}. In addition, we introduce a new objective function for $k$-BMF under which the problem becomes computationally easier and we explore the relationship between this new objective function and the original squared Frobenius distance.
Finally, we demonstrate that our proposed solution method is able to prove optimality for smaller datasets, while for larger datasets it provides solutions with better accuracy than the state-of-the-art heuristic methods. 
In addition, the entry-wise modelling of $k$-BMF in our formulations naturally extends to handle matrices with missing entries and perform binary matrix completion, we illustrate this way of application experimentally.

The rest of this paper is organised as follows. 
In Section \ref{section_formulations} we detail the three IP formulations for $k$-BMF and prove several results about their LP-relaxations. In Section \ref{section_new_objective}, we introduce a new objective function and explore its relation to the original squared Frobenius objective.
In Section \ref{section_column_generation} we detail a framework based on the large scale optimisation technique of column generation for the solution of our exponential formulation and discuss heuristics for the arising pricing problems. 
Finally, in Section \ref{section_experiments} we demonstrate the practical applicability of our approach on several artificial and real world datasets.

\section{Formulations.}
\label{section_formulations}

Given a binary matrix $X\in \B^{n\times m}$ and a fixed positive integer $k\ll \min(n,m)$ we wish to find two binary matrices $A\in \B^{n\times k}$ and $B\in \B^{k \times m}$
so that $\|X-Z \|^2_F$ is minimised, where $Z$ is the product of $A$ and $B$ and $\|\cdot\|_F$ denotes the Frobenius norm. Let $E=\{ (i,j) : x_{ij}=1\} \subset [n] \times [m]$ denote the index set of nonzero entries of $X$ where $[n] : = \{ 1, \dots, n \}$.
Both $X$ and $Z$ being binary matrices, the squared Frobenius and the entry-wise $\ell_1$ norm coincide and we may expand the objective function to get a linear expression
\begin{equation}
	\label{frob_obj}
	\|X-Z\|^2_F = \sum_{i=1}^n \sum_{j=1}^m |x_{ij} - z_{ij}|= \sum_{(i,j) \in E} (1-z_{ij}) +\sum_{(i,j) \not \in E} z_{ij}.
\end{equation}
For an incomplete binary matrix $X$ with missing entries, the above objective is slightly changed to $\sum_{(i,j) \in E} (1-z_{ij}) +\sum_{(i,j) \in \overline{E}}z_{ij}$ where $\overline{E} = \{(i,j): x_{ij}=0  \}$, to emphasise that $E\cup \overline{E}\not = [n]\times [m]$, and the factorisation error is only measured over known entries.
In the following sections we present three different integer programs for $k$-BMF all with the above derived linear objective function. 

\subsection{Compact formulation.}\label{section_compact_formulation}
We start with a formulation that uses a polynomial number of variables and constraints where we denote the McCormick envelope \cite{McCormick:1976} of $a,b\in[0,1]$ by
\begin{equation}
	MC(a,b) = \{y \in \R : \, 0\le y,\,a + b - 1 \le y, \, y\le a,\, y\le b \}.
\end{equation}
Note that if $a,b\in\B$ then  $MC(a,b)$ only contains the point $ab\in \B$ corresponding to the product of $a$ and $b$.
The following Compact Integer linear Program (CIP) models the entries of matrices $A,B,Z$ directly via binary variables $a_{i\ell}$, $b_{\ell j}$ and $z_{ij}$ respectively (for $i\in[n],\ell\in[k], j\in[m]$) and uses McCormick envelopes to avoid the appearance of quadratic terms that would correspond to the constraints $y_{i\ell j} = a_{i\ell} b_{\ell j}$,
\begin{align}
	(\text{CIP}) \quad \zeta_{\text{CIP}}= \min_{a,b,y,z} \;& \sum_{(i,j) \in E} (1-z_{ij}) +\sum_{(i,j)  \in \overline{E}} z_{ij}  \label{compact_obj}\\
	\text{s.t. }
	& y_{i\ell j}\le z_{ij} \le \sum_{l=1}^k y_{i l j} & i\in[n], j\in[m], \ell\in[k],\label{compact_boolean}\\
	& y_{i\ell j}\in MC(a_{i\ell},b_{\ell j})& i\in[n], j\in[m],\ell\in[k], \label{compact_McCormick}\\
	& a_{i\ell},b_{\ell j},z_{ij} \in \B &i\in[n], j\in[m],\ell\in[k]. \label{compact_integrality}
\end{align}
Constraints \eqref{compact_boolean} encode Boolean matrix multiplication, while a simple modification of the model in which constraints \eqref{compact_boolean} are replaced by $ z_{ij} = \sum_{\ell=1}^k y_{i \ell j}$ models $k$-BMF under standard arithmetic. The McCormick envelopes in constraints \eqref{compact_McCormick} ensure that for $a_{i \ell},b_{\ell j}\in \B$, $y_{i\ell j}$ are binary variables taking the value $a_{i\ell} b_{\ell j}$. Due to the objective function, constraints \eqref{compact_boolean} and the binary nature of $y_{i\ell j}$, the binary constraints on variables $z_{ij}$ may be relaxed to $z_{ij}\in[0,1]$ without altering optimal solutions of the formulation. 

The LP relaxation of CIP (CLP) is obtained by replacing constraints \eqref{compact_integrality}  by $a_{i\ell}, b_{\ell,j}, z_{ij} \in [0,1]$.
For $k=1$, we have $z_{ij}=y_{i1j}$ and 
the feasible region of CIP is the  Boolean Quadric Polytope (BQP) over a bipartite graph \cite{Padberg:1989:BQP:70486.70495}. 
The LP relaxation of BQP has half-integral vertices \cite{Padberg:1989:BQP:70486.70495}, which implies that  CLP for $k=1$ has half-integral vertices as well. One can show that in this case, a simple rounding in which fractional values of CLP are rounded down to $0$ gives a $2$-approximation to $1$-BMF \cite{Shi:2014:6899417}. This however, does not apply for $k>1$.
We next show that CLP for $k>1$ has an objective function value $0$.
\begin{proposition}
	\label{lemma_compact_fromulation_lp_0_obj} Given a binary matrix $X\in \B^{n\times m}$,  \textup{CLP} has optimal objective value $0$ for $k>1$. 
	Moreover, for $k>2$ \textup{CLP} has at least $k |E|+1$ 
	vertices with objective value $0$.
\end{proposition}
\proof{Proof.} 
	For each $(i,j)\in E$ let $L_{(i,j)}\subseteq [k]$ such that $|L_{(i,j)}|\ge 2$ and consider the point
	\begin{align*}
		a_{i\ell} &= \frac{1}{2}  \quad  i\in[n],\ell \in[k],	&b_{\ell j} &= \frac{1}{2} \quad \ell \in[k],  j\in[m],\\[.3cm]
		y_{i\ell j} &= \begin{cases}\frac{1}{2}&(i,j)\in E, \ell \in L_{(i,j)}\\
			0&\text{otherwise},\end{cases}
		&z_{ij} &= \begin{cases}1&(i,j)\in E, \\0&\text{otherwise}.\end{cases}
	\end{align*}
	For all $(i,j)\in[n]\times [m]$ and $\ell\in[k]$,  setting $a_{i\ell}=b_{\ell j}=\frac{1}{2}$ implies that $y_{i\ell j}\in  MC(\frac{1}{2}, \frac{1}{2}) =[0,\frac{1}{2}]$
	and  $ \sum_{l=1}^k y_{i l j}\ge 1$ holds for all $(i,j)\in E$, hence this point gives a feasible solution to  CLP with objective value $0$.
	For $k=2$, we can only set $L_{(i,j)}=[2]$ for all $(i,j)\in E$, hence the above construction leads to a single unique point.
	For $k>2$ however, as the choice of $L_{(i,j)}$'s is arbitrary, there are many feasible points with objective value $0$ of this form.  
	As each of these points can differ at only $k \, |E|$ entries corresponding to  entries $y_{i\ell j}$ for $(i,j)\in E$, $\ell\in [k]$, there are at most $k\, |E|+1$ affinely independent points among them. 
	Next we present $k \,|E|+1$ affinely independent points of this form. Since the objective value is $0$ at these points, they must lie on a face of dimension at least $k \,|E|$ and this face must have at least $k \,|E|+1$ vertices of CLP with objective value $0$.
	For each  $(i,j)^*\in E$ and $\ell^* \in [k]$, letting $L_{(i,j)}=[k]$ for all $ (i,j)\in E\setminus \{(i,j)^*\}$ and  $L_{(i,j)^*}=[k]\setminus\{\ell^* \}$ provides $k \,|E|$ different points of the above form. Each such point has exactly one entry $y_{i\ell j}$ along the indices $(i,j)\in E, \ell \in [k]$ which is zero.
	Hence the matrix whose columns correspond to these $k \,|E|$ points has a square submatrix of the form $\frac{1}{2}(J_{k |E|} - I_{k |E|})$ corresponding to entries $y_{i\ell j}$ for $(i,j)\in E, \ell \in [k]$, where $J_t$ is the all ones matrix of size $t\times t$ and $I_t$ is the identity matrix of size $t$. Since matrix $J_t - I_t$ is nonsingular, the $k \,|E|$ points are linearly independent. 
	In addition, letting $L_{(i,j)}=[k]$ for all $(i,j)\in E$ gives an additional point for which $y_{i\ell j}=\frac{1}{2}$ for all $(i,j)\in E, \ell \in [k]$, hence the corresponding part of this point is $\frac{1}{2}\mvec{1}$.
	Now subtracting $\frac{1}{2}\mvec{1}$ from the columns of $\frac{1}{2}(J_{k |E|} - I_{k |E|})$, we get the nonsingular matrix $-\frac{1}{2} I_{k |E|}$, hence the $k\, |E| +1 $ above constructed points are affinely independent.
	\Halmos
\endproof
The above result suggests that unless the factorisation error is $0$ i.e. the input matrix is of Boolean rank less than or equal to $k$, before improving the LP bound of CIP many fractional vertices need to be cut off.
Furthermore, for $k>1$, any feasible rank-$k$ factorisation $A\circ B$ and a permutation matrix $P\in \B^{k\times k}$ provide another feasible solution $AP\circ P^\top B$ with the same objective value. Hence, 
CIP is highly symmetric for $k>1$. These properties of CIP make it unlikely to be solved to optimality for $k>1$ in a reasonable amount of time for a large matrix $X$, though some symmetries 
may be broken by enforcing lexicographic ordering of rows of $B$.
For small matrices however, CIP constitutes the first approach to get optimal solutions to $k$-BMF.

\subsection{Exponential formulation I.}
Any $n\times m$ Boolean rank-$k$ matrix can be equivalently written as the Boolean combination of $k$ rank-$1$ binary matrices $\bigvee_{\ell=1}^k \mvec{a}_\ell \mvec{b}_\ell^\top$ for some $\mvec{a}_\ell\in\B^n,\mvec{b}_\ell\in\B^m$. This suggest to directly look for $k$ rank-1 binary matrices instead of introducing variables for all entries of factor matrices $A$ and $B$. The second integer program we detail for $k$-BMF relies on this approach by considering an implicit enumeration of rank-$1$ binary matrices.
Let $\RR$ denote the set of all rank-$1$ binary matrices of dimension $n\times m$ and let $\RR_{(i,j)}$ denote the subset of rank-$1$ matrices of $\RR$ which have the $(i,j)$-th entry equal to $1$,
\begin{align}
	\RR &:= \{\mvec{a} \mvec{b}^\top: \mvec{a}\in \B^n, \mvec{b}\in \B^m , \mvec{a},\mvec{b} \neq \mvec{0}\} \subset \B^{n\times m},\\
	\RR_{(i,j)} &:= \{ \mvec{a}\mvec{b}^\top\in \RR : a_i=b_j=1\} \subset \RR \;\;\quad\quad \qquad i\in[n],j\in[m].
\end{align}
Introducing a binary variable $q_r$ for each rank-$1$ matrix $r$ in $\RR$ and variables $z_{ij}$ corresponding to the  known entries of the $X$, we obtain the following Master Integer linear Program (MIP),
\begin{align}
	(\text{MIP}_{\text{F}}) \quad \zeta_{\text{MIP}}= \min_{z,q} \;& \sum_{(i,j) \in E} (1-z_{ij}) +\sum_{(i,j) \in \overline{E}} z_{ij}   \label{MIP_objective}\\
	\text{s.t. }
	& z_{ij} \le \sum_{ r\in \RR_{(i,j)} } q_r &  (i,j)\in E \label{MIP_upper_bool}\\
	&\sum_{ r\in \RR_{(i,j)} } q_r \le k\, z_{ij} &  (i,j) \in \overline{E}\label{MIP_lower_bool}\\
	& \sum_{r\in \RR} q_r \le k &\label{MIP_cardinality}\\
	& z_{ij}, q_r \in \B & (i,j)\in E\cup \overline{E},\; r\in\RR \label{MIP_bounds}
\end{align}
The objective, as before, measures the factorisation error in squared Frobenius norm, and subscript F in $\text{MIP}_{\text{F}}$ stands for Frobenius. Constraints \eqref{MIP_upper_bool} and \eqref{MIP_lower_bool} enforce Boolean matrix multiplication: $z_{ij}$ takes value $1$ if there is at least one active rank-1 binary matrix that covers entry $(i,j)$, otherwise it takes value $0$. Notice, that due to the difference in sign of objective coefficients for variables $z_{ij}$ with $(i,j)\in E$ and $(i,j)\in \overline{E}$ it is enough to declare constraints \eqref{MIP_upper_bool} and \eqref{MIP_lower_bool} for indices $(i,j)\in E$ and $(i,j) \in \overline{E}$ respectively.  Constraint \eqref{MIP_cardinality} ensures that at most $k$ rank-1 binary matrices are active and hence we get a rank-$k$ factorisation of $X$. 
Observe that constraints \eqref{MIP_upper_bool} together with $q_r$ being binary imply that $z_{ij}$ automatically takes binary values for $(i,j)\in E$, and due to the objective function it always takes the value at its upper bound, hence $z_{ij}\in\B$ may be replaced by $z_{ij}\le 1$ for all $(i,j)\in E$ without altering the optimum.
In contrast, $z_{ij}$ for $(i,j)\in \overline{E}$ need to be explicitly declared binary as 
otherwise, if there are some active rank-$1$ matrices ($q_r>0$) which cover a zero of $X$ ($r\in\RR_{(i,j)}$, $(i,j)\in \overline{E}$) then variable $z_{ij}$ corresponding to that zero takes the possibly fractional value $\frac{1}{k} \sum_{r\in \RR_{(i,j)}}q_r$.
One can also consider a \textit{strong formulation} of $\text{MIP}_{\text{F}}$ with exponentially many constraints, in which constraints \eqref{MIP_lower_bool} are replaced by $q_r\le z_{ij}$ for all $r\in \RR_{(i,j)}$ and $(i,j)\in\overline{E} $.

The LP relaxation of $\text{MIP}_{\text{F}}$ ($\text{MLP}_{\text{F}}$) is obtained by replacing the integrality constraints by $ z_{ij}, q_r \in [0,1]$.
Unlike CLP, the optimal objective value of $\text{MLP}_{\text{F}}$ ($\zeta_{\text{MLP}}$) is not always zero.
By comparing the rank of the factorisation, $k$ to the \textit{isolation number} of the input matrix $X$ we can deduce when $\text{MLP}_{\text{F}}$ will take non-zero objective value. 
We next give an extension of the definition of isolation number for binary matrices presented in \cite[Section 2.3]{Monson:1995}. 
\begin{definition}
	\label{isol_def}
	Let $X$ be a binary matrix with possibly missing entries. A set $S\subseteq E=\{ (i,j):x_{ij}=1  \}$ is said to be an \textit{isolated set of ones} if whenever $(i_1,j_1),(i_2,j_2)$ are two distinct members of $S$ then 
	\begin{enumerate*}[label=(\alph*)]
		\item $i_1\not = i_2$, $j_1 \not =j_2$ and
		\item $(i_1,j_2) \in \overline{E}$ or $(i_2,j_1) \in \overline{E}$ or both.
	\end{enumerate*}
	The size of the largest cardinality isolated set of ones of $X$ is denoted by $i(X)$ and is called the \textit{isolation number} of $X$.
\end{definition}
From the definition it follows that members of an isolated set of ones cannot be covered by a common rank-1 submatrix, and hence the isolation number provides a lower bound on the Boolean rank. 
The following result shows that $\text{MLP}_{\text{F}}$ must have non-zero objective value whenever $k$, the rank of the factorisation, is chosen so that it is strictly smaller than the isolation number.
\begin{proposition}
	\label{lemma_MLP_non_zero}
	Let $X$ have isolation number $i(X)>k$, then $\zeta_{\textup{MLP}}\ge \frac{1}{k} \left( i(X) - k \right)$.
\end{proposition}
\proof{Proof.}
	Let $S$ be an isolated set of ones of $X$ of cardinality $i(X)$. We will establish a feasible solution to the dual of  $\text{MLP}_{\text{F}}$ ($\text{MDP}_{\text{F}}$) with objective value $\frac{1}{k} \left( i(X)-k \right)$ implying the result. 
	
	Let us apply a change of variables $\xi_{ij}=1-z_{ij}$ for $(i,j)\in E$ for the ease of avoiding the constant term in the objective function of $\text{MLP}_{\text{F}}$. Then the bound constraints of $\text{MLP}_{\text{F}}$ can be written as $\xi_{ij} \ge 0$ for $(i,j)\in E$, $z_{ij}\ge 0$ for $(i,j)\in \overline{E}$ and $q_r\ge 0$, $r\in \RR$ as the objective function is minimising both $\xi_{ij}$ and $z_{ij}$ and we have the cardinality constrains on $q_r$.
	Associating dual variables $p_{ij}\ge 0$ $(i,j)\in E$ with constraints $\sum_{r\in \RR_{i,j}} q_r + \xi_{ij} \ge 1$, $s_{ij} \ge 0$ $(i,j) \in \overline{E}$ with constraints \eqref{MIP_lower_bool} and $\mu\ge0$ with constraint \eqref{MIP_cardinality}, the Master Dual Program ($\text{MDP}_{\text{F}}$) of $\text{MLP}_{\text{F}}$ is
	\begin{align}
		(\text{MDP}_{\text{F}})\;\;\zeta_{\text{MDP}} =\max_{p,s,\mu} \;& \sum_{(i,j)\in E}p_{ij} - k \,\mu\\
		\text{s.t. }
		&\sum_{(i,j)\in E\cap \supp(R)  } p_{ij} -  \sum_{(i,j) \in \overline{E} \cap \supp(R) } s_{ij} \le\mu  & R\in \RR \label{MDP_balance_con}\\
		& 0 \le p_{ij} \le 1& (i,j)\in E\\
		& 0 \le s_{ij} \le \frac{1}{k}& (i,j)\in \overline{E}\\
		&0\le\mu,
	\end{align}
	where $\supp (R) = \{(i,j) : r_{ij}=1  \}$.
	
	Let $s_{ij}=\frac{1}{k}$ for $(i,j)\in \overline{E}$ and  let $p_{ij}=\frac{1}{k}$ for $(i,j)\in S$ and $p_{ij}=0$ for all other $(i,j)\in E\setminus S$. 
	The bound constraints on $p_{ij}$ and $s_{ij}$ are satisfied then. 
	It remains to choose $\mu\ge 0$ such that we satisfy constraint \eqref{MDP_balance_con}
	for all rank-$1$ binary matrices $R\in \RR$. 
	Let $R\in \RR$ be a submatrix of $X$, so we have $| \overline{E}\cap \supp(R) |=0$. Then by the definition of isolated sets of ones, $R$ can contain at most one element from $S$ and hence we have $|\supp(R) \cap S| \le 1$. This tells us that for any $\mu\ge \frac{1}{k}$, constraint \eqref{MDP_balance_con} is satisfied for all $R\in\RR$ that is a submatrix of $X$.
	Now let $R\in\RR$ be a rank-$1$ binary matrix which covers at least one zero entry of $X$. Then $R$ may contain more than one element from $S$. However, if it contains more than one element from $S$  then it must also contain at least $\binom{|\supp(R)\cap S|}{2}$-many zeros as for any two distinct elements $(i_1,j_1),(i_2,j_2)$ in $S$ we have $(i_1,j_2) \in \overline{E}$ or $(i_2,j_1) \in \overline{E}$ by the definition of isolated set of ones. Hence, for all $R\in \RR$ such that $|\overline{E} \cap \supp(R)|>0$, constraint \eqref{MDP_balance_con}  satisfies 
	\begin{equation}
		\frac{1}{k} | S\cap \supp(R)| - \frac{1}{k} |\overline{E} \cap \supp(R)| \le\frac{1}{k} | S\cap \supp(R)| -  \frac{1}{k}\binom{| S\cap \supp(R)|}{2} \le \frac{1}{k}.
	\end{equation}
	Thus we can set $\mu =\frac{1}{k}$ to get the objective value $\frac{1}{k}\left( i(X) - k \right) \le \zeta_{\text{MDP}} =\zeta_{\text{MLP}}$, which provides a non-zero bound on $\text{MLP}_{\text{F}}$ for all $k< i(X)$.
	\Halmos
\endproof
The following example shows that we cannot strengthen Proposition \ref{lemma_MLP_non_zero} by replacing the condition $k<i(X)$  with the requirement that $k$ has to be strictly smaller than the Boolean rank of $X$.
\begin{example}
	\label{example_comp_I_4}
	Let $X=J_4 - I_4$, where $J_4$ is the $4 \times 4$ matrix of all $1$s and $I_4$ is the $4\times 4$ identity matrix. One can verify that the Boolean rank of $X$ is $4$ and its isolation number is $3$. 
	For $k=3$, the optimal objective value of $\text{MLP}_{\text{F}}$ is $0$ which is attained by a fractional solution in which the following $6$ rank-$1$ binary matrices are active with weight $\frac12$.
	
	\setlength{\tabcolsep}{12pt} 
	\begin{tabular}{cccccc}
		$q_1=\frac{1}{2}$ & $q_2=\frac{1}{2}$ & $q_3=\frac{1}{2}$ & $q_4=\frac{1}{2}$ &$q_5=\frac{1}{2}$ & $q_6=\frac{1}{2}$ \\
		$\begin{bmatrix}
		0 & 0 & 0 & 0  \\
		1 & 0 & 1 & 0  \\
		0 & 0 & 0 & 0  \\
		1 & 0 & 1 & 0  
		\end{bmatrix}$
		&
		$\begin{bmatrix}
		0 & 1 & 1 & 0  \\
		0 & 0 & 0 & 0  \\
		0 & 0 & 0 & 0  \\
		0 & 1 & 1 & 0  
		\end{bmatrix}$
		&
		$\begin{bmatrix}
		0 & 1 & 0 & 1  \\
		0 & 0 & 0 & 0  \\
		0 & 1 & 0 & 1  \\
		0 & 0 & 0 & 0  
		\end{bmatrix}$
		&
		$\begin{bmatrix}
		0 & 0 & 0 & 0  \\
		1 & 0 & 0 & 1  \\
		1 & 0 & 0 & 1  \\
		0 & 0 & 0 & 0  
		\end{bmatrix}$
		&
		$\begin{bmatrix}
		0 & 0 & 1 & 1  \\
		0 & 0 & 1 & 1  \\
		0 & 0 & 0 & 0  \\
		0 & 0 & 0 & 0  
		\end{bmatrix}$
		&
		$\begin{bmatrix}
		0 & 0 & 0 & 0  \\
		0 & 0 & 0 & 0  \\
		1 & 1 & 0 & 0  \\
		1 & 1 & 0 & 0  
		\end{bmatrix}$
	\end{tabular}
\end{example}

\subsection{Exponential formulation II.}

For $t\in [2^m-1]$ let $\mvec{\beta}_t\in\B^{m}$ be the vector denoting the binary encoding of $t$ and note that these vectors give a complete enumeration of all non-zero binary vectors of size $m$.
Let $\beta_{tj}$ denote the $j$-th entry of $\mvec{\beta}_t$.
In \cite{Lu:2008:OBM:1546682.1547186}, the authors present the following Exponential size Integer linear Program (EIP) formulation using a separate indicator variable $d_t$ for each one of these exponentially many binary vectors $\mvec{\beta}_t$,
\begin{align}
	(\text{EIP}) \quad \zeta_{\text{EIP}}= \min_{\alpha,z,d} \;& \sum_{(i,j) \in E} (1-z_{ij}) +\sum_{(i,j)  \in \overline{E}} z_{ij}  \\
	\text{s.t. }
	& z_{ij} \le \sum_{t=1}^{2^m-1} \alpha_{i t} \beta_{t j}& (i,j)\in E, \label{exp_bool_upper}\\
	&\sum_{t=1}^{2^m-1} \alpha_{i t} \beta_{t j}\le k z_{ij} & (i,j) \in \overline{E}, \label{exp_bool_lower}\\
	& \sum_{t=1}^{2^m-1} d_{t}\le k & \label{exp_k}\\
	& \alpha_{i t}\le d_{t} & i\in[n], t\in[2^m-1], \label{exp_alpha}\\
	& z_{ij}, d_t, \alpha_{i t}\in \B & (i,j)\in E \cup \overline{E}, t\in[2^m-1]. \label{exp_int}
\end{align}
The above formulation has an exponential number of variables and constraints but it is an  integer  linear program as $\beta_{tj}$ are input parameters to the model. 
Let ELP be the LP relaxation of EIP. Observe that due to the objective function the bound constraints in ELP may be simplified to $z_{ij},\alpha_{it}, d_t\ge 0$  for all $i,j,t$ and $z_{ij}\le1 $ for $(i,j)\in E$ without changing the optimum.
To solve EIP or ELP explicitly, one needs to enumerate all binary vectors $\mvec{\beta}_t$, which is possible only up to a very limited size. 
To the best of our knowledge, no method is available that avoids explicit enumeration and can guarantee the optimal solution of EIP. Previous attempts at computing a rank-$k$ factorisation via EIP all relied on working with only a small heuristically chosen subset of vectors $\mvec{\beta}_t$ \cite{Lu:2008:OBM:1546682.1547186,Lu:2014:9400730720140101}. However, if there was an efficient method to solve ELP, the following result shows it to be as strong as the LP relaxation of $\text{MIP}_{\text{F}}$.
\begin{proposition}
	\label{lemma_LPexp_equiv_to_MLP_F}
	The optimal objective values of \textup{ELP} and $\textup{MLP}_{\textup{F}}$ are equal.
\end{proposition}
\proof{Proof.}
	Note that due to constraints \eqref{MIP_upper_bool} and \eqref{MIP_lower_bool} in $\text{MLP}_{\text{F}}$ and  constraints \eqref{exp_bool_upper} and \eqref{exp_bool_lower} in ELP,  it suffices to show that for any feasible solution $\alpha_{i t},d_t$  of ELP one can build a feasible solution $q_r$ of $\text{MLP}_{\text{F}}$ for which $\sum_{t=1}^{2^m-1}\alpha_{it} \beta_{tj} = \sum_{r\in\RR_{(i,j)}}q_r$, and vice-versa.
	
	First consider a feasible solution $\mvec{\alpha}_t\in\R^n,~d_t\in\R$ (for $t\in[2^m-1]$)  to ELP and note that by  constraint \eqref{exp_alpha} 
	we have $0\le \alpha_{it}\le d_t$ for all $i\in[n]$ and $t\in [2^m-1]$.
	We can therefore express each $\mvec{\alpha}_t$ as a convex combination of binary vectors in $\B^{n}$ scaled by $d_t$,
	\begin{equation}
		\mvec{\alpha}_t = d_t \,\sum_{s=1}^{2^n-1} \lambda_{s,t} \; \mvec{a}_s 
		\quad \mvec{a}_s\in\B^{n}\setminus\{\mvec{0} \}, \quad 
		\sum_{s=1}^{2^n-1}\lambda_{s,t} \le 1, \quad \lambda_{s,t} \ge 0, \quad s\in [2^n-1]
	\end{equation}
	where $\mvec{a}_s$ denotes the binary encoding of $s$. Note that we do not require $\lambda_{s,t}$'s to add up to 1 as we exclude the zero vector.  
	We can therefore rewrite the solution of ELP as follows
	\begin{equation}
		\sum_{t=1}^{2^m-1} \mvec{\alpha}_t \mvec{\beta}^\top_t  =  \sum_{t=1}^{2^m-1} \sum_{s=1}^{2^n-1}  d_t \,\lambda_{s,t} \,\mvec{a}_s \mvec{\beta}_t^\top = \sum_{t=1}^{2^m-1} \sum_{s=1}^{2^n-1} q_{s,t}\mvec{a}_s \mvec{\beta}_t^\top 		\quad \text{where } q_{s,t} := d_t\,\lambda_{s,t}.
	\end{equation}
	Now it is easy to see that $\mvec{a}_s\mvec{\beta}_t^\top\in \RR$ and since $\sum_{t=1}^{2^m-1} d_t \le k$ holds in any feasible solution to ELP, we get $\sum_{s=1}^{2^n-1} \sum_{t=1}^{2^m-1}  q_{s,t} \le k$, which shows that $q_{s,t}$ is feasible for $\text{MLP}_{\text{F}}$. 
	
	The construction works backwards as well, as any feasible solution to $\text{MLP}_{\text{F}}$ 
	can be written as \\
	$\sum_{s=1}^{2^n-1}\sum_{t=1}^{2^m-1} q_{s,t} \mvec{a}_s \mvec{\beta}_t^\top$ for some rank-$1$ binary matrices $\mvec{a}_s \mvec{\beta}_t^\top\in \RR$ and corresponding variables $q_{s,t}\ge 0$. Now let $\mvec{\alpha}_t := \sum_{s=1}^{2^n-1} q_{s,t}\, \mvec{a}_s$ and $d_t := \max_{i\in [n]} \alpha_{it} $ to satisfy $\alpha_{it}\le d_t$. Then since we started from a feasible solution to $\text{MLP}_{\text{F}}$, we have $\sum_{s=1}^{2^n-1}\sum_{t=1}^{2^m-1} q_{s,t}\le  k$ and hence $\sum_{t=1}^{2^m-1} d_t \le k$ is satisfied too.
	\Halmos
\endproof


\section{Working under a new objective}
\label{section_new_objective}

\begin{figure}[!ht]
		\caption{\label{fig_color_rectangles}Example matrices for which $\zeta(1)=k\, \zeta_F$}
		\centering
		\begin{subfigure}[b]{.2\textwidth}
			\centering
			\begin{tikzpicture}[scale=0.3]
			\pgfmathsetmacro\y{5}
			\pgfmathsetmacro\m{.5}
			\pgfmathsetmacro\x{2}
			\filldraw[gray!20!white, draw=black,thick] (-\y,\y) rectangle (\m,-\m);
			\filldraw[gray!40!white, draw=black,thick] (\y,-\y) rectangle (-\m,\m);
			\filldraw[white!40!white, draw=black,thick] (\m,-\m) rectangle (-\m,\m);
			\end{tikzpicture}
			\caption{$k=2$\label{example_k_2}}
		\end{subfigure}
		\begin{subfigure}[b]{.3\textwidth}
			\centering
			\begin{tikzpicture}[scale=0.3]
			\pgfmathsetmacro\y{4.5}
			\pgfmathsetmacro\m{.5}
			\pgfmathsetmacro\x{2.5}
			%
			\filldraw[gray!20!white, draw=black,thick] (-\y,\y) rectangle (\m,-\m);
			\filldraw[gray!40!white, draw=black,thick] (\y,-\y) rectangle (-\m,\m);
			\filldraw[blue!40!white, draw=black,thick] (-\m+\x,\m-\x) rectangle (-\m,\m);
			\filldraw[blue!40!white, draw=black,thick] (-\y-\x,-\x+\m) rectangle (-\y,\m);
			\filldraw[blue!40!white, draw=black,thick] (-\m,\y) rectangle (-\m+\x,\y+\x);
			\filldraw[blue!40!white, draw=black,thick] (-\y,\y) rectangle (-\x-\y,\y+\x);
			\filldraw[green!40!white, draw=black,thick] (\m-\x,-\m+\x) rectangle (\m,-\m);
			\filldraw[green!40!white, draw=black,thick] (\y+\x,\x-\m) rectangle (\y,-\m);
			\filldraw[green!40!white, draw=black,thick] (\m,-\y) rectangle (\m-\x,-\y-\x);
			\filldraw[green!40!white, draw=black,thick] (\y,-\y) rectangle (\x+\y,-\y-\x);
			%
			\filldraw[white!40!white, draw=black,thick] (\m,-\m) rectangle (-\m,\m);
			\draw (\y,-\y) rectangle (\y-\x,-\y+\x);
			\draw (-\y,\y) rectangle (-\y+\x,\y-\x);
			\end{tikzpicture}
			\caption{$k=4$}
			\label{example_k_4}
		\end{subfigure}
		\begin{subfigure}[b]{.4\textwidth}
			\centering
			\begin{tikzpicture}[scale=0.2]
			\pgfmathsetmacro\y{8}
			\pgfmathsetmacro\m{.5}
			\pgfmathsetmacro\x{3}
			%
			\filldraw[gray!20!white, draw=black,thick] (-\y,\y) rectangle (\m,-\m);
			\filldraw[gray!40!white, draw=black,thick] (\y,-\y) rectangle (-\m,\m);
			\filldraw[blue!40!white, draw=black,thick] (-\m+\x,\m-\x) rectangle (-\m,\m);
			\filldraw[blue!40!white, draw=black,thick] (-\y-\x,-\x+\m) rectangle (-\y,\m);
			\filldraw[blue!40!white, draw=black,thick] (-\m,\y) rectangle (-\m+\x,\y+\x);
			\filldraw[blue!40!white, draw=black,thick] (-\y,\y) rectangle (-\x-\y,\y+\x);
			\filldraw[green!40!white, draw=black,thick] (\m-\x,-\m+\x) rectangle (\m,-\m);
			\filldraw[green!40!white, draw=black,thick] (\y+\x,\x-\m) rectangle (\y,-\m);
			\filldraw[green!40!white, draw=black,thick] (\m,-\y) rectangle (\m-\x,-\y-\x);
			\filldraw[green!40!white, draw=black,thick] (\y,-\y) rectangle (\x+\y,-\y-\x);
			%
			\filldraw[pink!100!white, draw=black,thick] (\x-\m,\m) rectangle (2*\x-2*\m,-\m); 
			\filldraw[pink!100!white, draw=black,thick] (2*\x-2*\m,-\x+\m) rectangle (\x-\m,-2*\x+2*\m); 
			\filldraw[pink!100!white, draw=black,thick] (-\m,-\x+\m) rectangle (\m,-2*\x+2*\m); 
			\filldraw[pink!50!white, draw=black,thick] (-\y-\x,-\m) rectangle (-\y-2*\x,\m); 
			\filldraw[pink!50!white, draw=black,thick] (-\y-\x, -2*\x+2*\m) rectangle (-\y-2*\x,-\x+\m); 
			\filldraw[pink!100!white, draw=black,thick] (-\m,\y+2*\x) rectangle (\m,\y+\x); 
			\filldraw[pink!100!white, draw=black,thick] (\x-\m,\y+2*\x) rectangle (2*\x-2*\m,\y+\x); 
			\filldraw[pink!100!white, draw=black,thick] (-\y-\x,\y+\x) rectangle (-2*\x-\y,\y+2*\x);
			%
			\filldraw[yellow!50!white, draw=black,thick] (-\x+\m,-\m) rectangle (-2*\x+2*\m,\m); 
			\filldraw[yellow!50!white, draw=black,thick] (-2*\x+2*\m,\x-\m) rectangle (-\x+\m,2*\x-2*\m); 
			\filldraw[yellow!50!white, draw=black,thick] (\m,\x-\m) rectangle (-\m,2*\x-2*\m); 
			\filldraw[yellow!100!white, draw=black,thick] (\y+\x,\m) rectangle (\y+2*\x,-\m); 
			\filldraw[yellow!100!white, draw=black,thick] (\y+\x, \x-\m) rectangle (\y+2*\x,2*\x-2*\m); 
			\filldraw[yellow!50!white, draw=black,thick] (\m,-\y-2*\x) rectangle (-\m,-\y-\x); 
			\filldraw[yellow!50!white, draw=black,thick] (-\x+\m,-\y-2*\x) rectangle (-2*\x+2*\m,-\y-\x); 
			\filldraw[yellow!50!white, draw=black,thick] (\y+\x,-\y-\x) rectangle (2*\x+\y,-\y-2*\x);
			%
			\filldraw[white!40!white, draw=black,thick] (\m,-\m) rectangle (-\m,\m);
			\draw (\y,-\y) rectangle (\y-\x,-\y+\x);
			\draw (-\y,\y) rectangle (-\y+\x,\y-\x);
			\end{tikzpicture}
			\caption{$k=6$}
			\label{example_k_6}
		\end{subfigure}
	\end{figure}
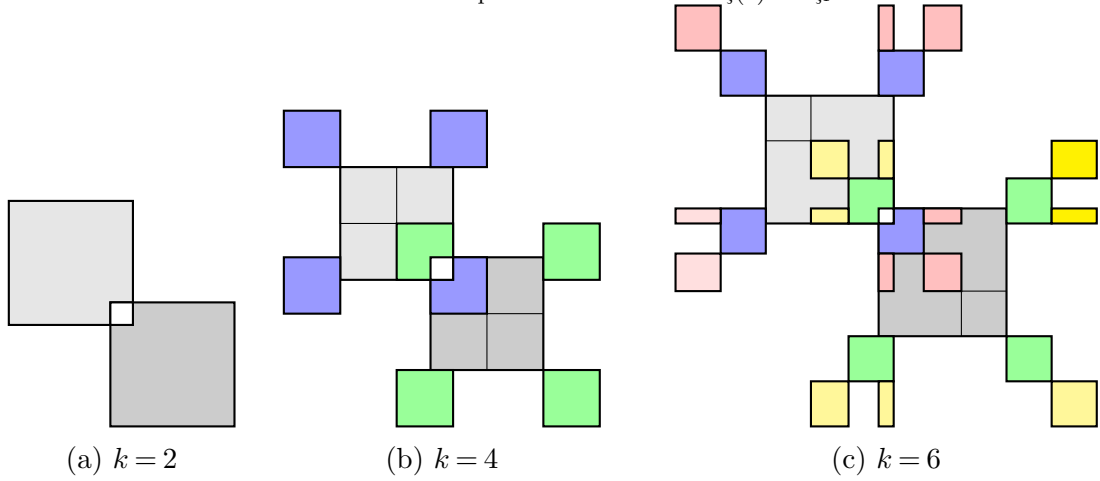

In the previous section, we presented formulations for $k$-BMF which measured the factorisation error in the squared Frobenius norm, which coincides with the entry-wise $\ell_1$ norm as showed in Equation \eqref{frob_obj}.
In this section, we explore another objective function which introduces an asymmetry between how false negatives and false positives are treated. 
Whenever a $0$ entry is erroneously covered in a rank-$k$ factorisation, it may be covered by up to $k$ rank-1 binary matrices. Our new objective function attributes an error term to each $0$ entry which is proportional to the number of rank-1 matrices covering that entry. As previously, by denoting $Z=A\circ B$ a rank-$k$ factorisation of $X$, the new objective function is
\begin{equation}
	\label{equation_new_objective}
	\zeta(\rho ) = \sum_{(i,j) \in E} (1-z_{ij} )+\rho \sum_{(i,j)  \in \overline{E}} \sum_{\ell=1}^k a_{i\ell}b_{\ell j}.
\end{equation}
Note that the  constraints $a_{i\ell} b_{\ell j}\le z_{ij} \le \sum_{\ell=1}^k a_{i\ell} b_{\ell j}$ encoding Boolean matrix multiplication imply that $\frac{1}{k}\sum_{\ell=1}^k a_{i\ell} b_{\ell j}\le z_{ij}  \le\sum_{\ell=1}^k a_{i\ell} b_{\ell j}$.
Therefore, denoting the original squared Frobenius norm objective function in Equation \eqref{frob_obj} by $\zeta_F$, for any $X$ and rank-$k$ factorisation $Z$ of $X$ the following relationship holds between $\zeta_{F}$ and $\zeta(1)$, $\zeta(\frac{1}{k})$,
\begin{align}
	\label{k_times}
	\zeta_F 
	&
	\le
	\zeta(1)
	\le
	\sum_{(i,j) \in E} (1-z_{ij}) + \sum_{(i,j)\in \overline{E}} k\,z_{ij}
	\le
	k \, \zeta_F
	\qquad\text{and}
	\qquad 
	\frac{1}{k} \zeta_F 
	\le
	\zeta(\frac{1}{k}) 
	\le 
	\zeta_F. 
\end{align}
We next show that this new objective function $\zeta(\rho)$ with $\rho=1$ can overestimate the original objective $\zeta_{F}$ by a factor of $k$. 
But first, we need a technical result which shows that whenever the input matrix $X$ contains repeated rows or columns we may assume that an optimal factorisation exists which has the same row-column repetition pattern. 
\begin{lemma}[Preprocessing]
	\label{lemma_preprocessing}
	Let $X$ contain some duplicate rows and columns. Then there exists an optimal rank-$k$ binary matrix factorisation of $X$ under objective $\zeta_F$ (or $\zeta(\rho)$) whose rows and columns corresponding to identical copies in $X$ are identical.
\end{lemma}
\proof{Proof.}
	Since the transpose of an optimal rank-$k$ factorisation is optimal for $X^\top$, it suffices to consider the rows of $X$.
	Furthermore, it suffices to consider only one set of repeated rows of $X$, 
	so let $I\subseteq [n]$ be the index set of a set of identical rows of $X$. We then need to show that there exists an optimal rank-$k$ factorisation whose rows indexed by $I$ are identical.
	Let $Z=A\circ B$ be an optimal rank-$k$ factorisation of $X$ under objective $\zeta_F$.
	For all $i_1,i_2\in I$ we must have
	\begin{equation}
		\label{eq_equal_contributiion_frob}
		\sum_{j: (i_1,j)\in E} (1-z_{ij}) + \sum_{j: (i_1,j)\in \overline{E}} z_{ij} = 
		\sum_{j: (i_2,j)\in E} (1-z_{ij}) + \sum_{j: (i_2,j)\in \overline{E}} z_{ij}
	\end{equation}
	as otherwise replacing $A_{i,:}$ for each $i\in I$ with row $A_{i^*,:}$ where $i^*\in I$ is a row index for which the above sum is minimised leads to a smaller error factorisation. Then since \eqref{eq_equal_contributiion_frob} holds, replacing $A_{i,:}$ for each $i\in I$ with row $A_{i^*,:}$ for any $i^*\in I$ leads to an optimal solution of the desired property. 
	Similarly, if $Z$ is an optimal factorisation under objective $\zeta(\rho)$, then for all $i_1,i_2\in I$ the corresponding objective terms must  equal and hence an optimal solution of the desired property exists. \Halmos
\endproof
This result implies that whenever the input matrix $X$ contains repeated rows or columns we may solve the following problem on a smaller matrix instead. 
Let $X' \in \B^{n' \times m'}$ be the binary matrix obtained from $X$ by replacing each duplicate row and column by a single representative and let $\mvec{r}\in \Z_+^{n'}$ and $\mvec{c}\in\Z_+^{m'}$ be the counts of each unique row and column of $X'$ in $X$ respectively.
Let $E'$ and $ \overline{E'}$ denote the non-zero and zero entry index sets of $X'$ respectively.
By Lemma \ref{lemma_preprocessing} an optimal rank-$k$ factorisation $Z'=A'\circ B'$ of $X'$ under the \textit{updated} objective function 
\begin{equation}
	\zeta_F':=\sum_{(i,j)\in E'} r_i \, c_j\, (1-z'_{ij}) + \sum_{(i,j) \in \overline{E'}} r_i \, c_j\, z'_{ij}
\end{equation}
(or  $\zeta'(\rho):=\sum_{(i,j)\in E'} r_i \, c_j\, (1-z'_{ij}) + \rho \sum_{(i,j) \in \overline{E'}} r_i \, c_j \sum_{\ell=1}^k\, a'_{i\ell} b'_{\ell j}$)
leads to an optimal rank-$k$ factorisation of $X$ under the original objective function $\zeta_{F}$ (or $\zeta(\rho)$).
\begin{proposition} \label{tight_example} For each positive integer $k$ there exists a matrix $X(k)$ for which the optimal rank-$k$ binary matrix factorisations under objectives $\zeta_F$ and $\zeta(1)$ satisfy $\zeta(1) = k \, \zeta_F$.
\end{proposition}
\proof{Proof.}
	The idea behind the proof is to consider a matrix $Z(k)$ of exact Boolean rank-$k$ in which all the $k$ rank-$1$ components (rectangles) overlap at a unique middle entry and then replace this entry with a $0$ to obtain $X(k)$. Now $X(k)$ and $Z(k)$ are exactly at distance $1$ in the squared Frobenius norm and hence $Z(k)$ is a rank-$k$ factorisation of $X(k)$ with objective value $\zeta_{F}=1$. On the other hand, since exactly $k$ rectangles cover the entry at which $X(k)$ and $Z(k)$ differ, if $Z(k)$ is taken as a rank-$k$ factorisation of $X(k)$ under objective $\zeta(1)$ it incurs an error of size $k$. Figure \ref{fig_color_rectangles} shows the idea how to build such a $X(k)$ for $k=2,4,6$. Each colour corresponds to a rank-1 component and white areas correspond to $0$s.
	We first consider the case when $k$ is even.
	For $k=2$ take the symmetric matrix $X(2)$ as in Equation \eqref{X_2} which corresponds to Figure \ref{example_k_2}. Since $X(2)$ has repeated rows and columns, according to Lemma \ref{lemma_preprocessing} we may simplify the problem by replacing $X(2)$ by $X'(2)$ and recording a weight vector for the rows and columns which indicate how many times each row and column is repeated. This weight vector is then used to update each entry in the objective function with the corresponding weight. For $X(2)$ the row and column weight vectors coincide as $X(2)$ is symmetric and we denote it by $\mvec{w}(2)$. 
	\begin{equation}
		\label{X_2}
		X(2)=
		{
			\begin{bmatrix}
				1 & 1 & 0\\
				1 & 1 & 0\\
				1 & 1 & 0\\
				1 & 0 & 1\\
				0 & 1 & 1\\
				0 & 1 & 1\\
				0 & 1 & 1
			\end{bmatrix}
			\circ 
			\begin{bmatrix}
				1 & 0 & 0\\
				1 & 0 & 0\\
				1 & 0 & 0\\
				0 & 1 & 0\\
				0 & 0 & 1\\
				0 & 0 & 1\\
				0 & 0 & 1
			\end{bmatrix}^\top
		}
		=
		{
			\begin{bmatrix}
				1 & 1 & 1 & 1 & 0 & 0 & 0 \\
				1 & 1 & 1 & 1 & 0 & 0 & 0 \\
				1 & 1 & \colorbox{gray!20}{1} &  1 & 0 & 0 & 0 \\
				1 & 1 &  1 & 0 & \colorbox{gray!20}{1} & 1 & 1 \\
				0 & 0 & 0 & \colorbox{gray!20}{1} &  1 & 1 & 1\\
				0 & 0 & 0 & 1 & 1 & 1 & 1\\
				0 & 0 & 0 & 1 & 1 & 1 & 1\\
			\end{bmatrix}
		}
		\Rightarrow
		X'(2)=
		{
			\begin{bmatrix}
				1 & 1 & 0  \\
				1 & 0 & 1  \\
				0 & 1 & 1 \\
			\end{bmatrix}\text{with }
			\mvec{w}(2)=
			\begin{bmatrix}
				3\\1\\3
		\end{bmatrix}}
	\end{equation}
	The Boolean rank of $X(2)$ is $3$, which one can confirm by looking at a size $3$ isolated set of ones (shadowed entries) and an exact rank-$3$ factorisation shown in Equation \eqref{X_2}. 
	Let $Z(2)$ be obtained from $X(2)$ by replacing the $0$ at entry $(4,4)$ by a $1$. 
	$Z(2)$ clearly has Boolean rank $2$, hence it is a feasible rank-$2$ factorisation of $X(2)$.
	Under objective $\zeta_F$ $Z(2)$ incurs an error of size $1$, which is optimal as $\zeta_F\ge 1$  by  $X(2)$ being of Boolean rank-$3$. 
	On the other hand, under objective $\zeta(1)$ $Z(2)$  has objective value $2$ as the middle entry is covered twice. To see that $Z(2)$ is optimal under $\zeta(1)$ observe that every entry in $X'(2)$ apart from the middle entry has weight strictly greater than $2$. Hence not covering a $1$ of $X'(2)$ or covering a $0$ different from the middle entry incurs an error strictly greater than $2$.

	For $k>2$ even let us give a recipe to construct a symmetric matrix $X'(k)$ and corresponding weight vector $\mvec{w}(k)$. Let $t = \frac{k}{2}-1$ and let the following $(4t + 3) \times (4t + 3)$ matrix be $X'(k)$, where $I_t$ is the identity matrix of size $t\times t$, $\tilde{I}_t$ is the reverted identity matrix of size $t\times t$ (so 
	$\tilde{I}_2=\left[
	\begin{smallmatrix} 
	0 & 1 \\
	1 & 0
	\end{smallmatrix} \right]$) and $J_t$ is the all ones matrix of size $t\times t$, 
	\begin{equation*}
		{
			X'(k)=
			\begin{bmatrix}	
				\colorbox{gray!20}{$I_{t}$} & &	& \mvec{1}_{t} 	&\tilde{I}_{t} &&\\
				& \colorbox{gray!20}{1}& \mvec{1}^\top_{t} & 1	&&&\\
				&\mvec{1}_{t}	& J_{t} & \mvec{1}_{t}     &&& \tilde{I}_{t}  \\
				\mvec{1}^\top_{t} &1 &\mvec{1}^\top_{t} & 0  &\mvec{1}^\top_{t}  &\colorbox{gray!20}{1} & \mvec{1}^\top_{t} \\
				\tilde{I}_{t}&&	&\mvec{1}_{t}     & J_{t} 	& 	\mvec{1}_{t} &				 \\
				&&	&\colorbox{gray!20}{1}   & 	\mvec{1}^\top_{t} &1& \\
				&	&	\tilde{I}_{t}  		&\mvec{1}_{t}  & & & \colorbox{gray!20}{$I_{t}$}	 \\ 
			\end{bmatrix}, \;\;
			\mvec{w}(k) =
			\begin{bmatrix}
				(k+1) \mvec{1}_t \\ (k+1) \\ (k+1)\, \mvec{1}_t \\ 1 \\  (k+1)\, \mvec{1}_t \\ (k+1) \\ (k+1) \mvec{1}_t
			\end{bmatrix},
		}\;\;
		{
			A'(k)=
			\begin{bmatrix}	
				I_{t} & &	&\\
				&1&&\\
				& \mvec{1}_{t} && \tilde{I}_{t}  \\
				\mvec{1}^\top_t & 1 & 1 & \mvec{1}^\top_t\\
				\tilde{I}_{t}&	&\mvec{1}_{t}   &				 \\
				&&1& \\
				&&&I_{t}		 
			\end{bmatrix}.
		}
	\end{equation*}
	$X'(k)$ has isolation number $i(X'(k))\ge 2t + 3 =k+1$ (indicated by the shadowed entries),  so no rank-$k$ factorisation can have zero error.
	Let $Z'(k)$ be obtained from $X'(k)$ by replacing the middle $0$ by a $1$ and let its weight vector be the same as of $X'(k)$. The Boolean rank of $Z'(k)$ is then at most $k$ as $Z'(k)=A'(k) \circ A'^\top (k)$ is an exact factorisation and $A'(k)$ is of dimension $(4t + 3) \times k$. This factorisation is illustrated in Figure \ref{fig_color_rectangles} for $k=4,6$. Therefore $Z'(k)$ is a feasible rank-$k$ factorisation of $X'(k)$.
	Now $Z'(k)$ under objective function $\zeta_{F}$ has error $1$ and hence it is optimal.
	In contrast, $Z'(k)$ evaluated under objective $\zeta(1)$ has error $k$ as the middle $0$ is covered $k$ times and it has weight $1$.
	To see that $Z'(k)$ is optimal under $\zeta(1)$ as well, note that all entries of $X'(k)$ apart from the middle $0$ have weight strictly greater than $k$. Therefore, any other rank-$k$ factorisation which does not cover a $1$ or covers a $0$ which is not the middle $0$, incurs an error strictly greater than $k$, and hence $Z'(k)$ is optimal under objective $\zeta(1)$ with value $k\cdot \zeta_F$.
	
	For $k=1$, all $1$-BMFs satisfy $\zeta_{F}=\zeta(1)$ by definition. For $k>1$ odd, we can obtain $X'(k)$ and $\mvec{w}(k)$ from $X'(k+1)$ and $\mvec{w}(k+1)$ by removing the first row and column of $X'(k+1)$ and the corresponding first entry of $\mvec{w}(k+1)$. For $X'(k)$ then, the same reasoning holds as for $k$ even.
	\Halmos
\endproof

While Proposition \ref{tight_example} shows that $\zeta(1)$ can be $k$ times larger than the Frobenius norm objective $\zeta_F$, the matrices in the proof are quite artificial, and in practice we observe that not many zeros are covered by more than a few rank-1 matrices.  
Therefore it is worth considering the previously introduced formulations for $k$-BMF with the new objective $\zeta(\rho)$.

Let us denote a modification of formulation $\text{MIP}_{\text{F}}$ with the new objective function $\zeta(\rho)$ as MIP($\rho$) and use the transformation $\xi_{ij} = 1-z_{ij}$ for $(i,j)\in E$ to get
\begin{align}
	(\text{MIP}(\rho))\;\;\zeta_{\text{MIP}(\rho)} = \min_{\xi,q}\; & \sum_{(i,j)\in E} \xi_{ij} +
	\rho \sum_{(i,j) \in \overline{E}} \sum_{r\in \RR_{(i,j)}} q_r \label{MIP_rho_objective}\\
	\text{s.t. }
	& \sum_{ r\in \RR_{(i,j)} } q_r + \xi_{ij} \ge 1 & (i,j)\in E \label{MIP_rho_upper_bool}\\
	&\sum_{ r\in \RR } q_r \le k & \label{MIP_rho_cardinality}\\
	&\xi_{ij} \ge 0,\; q_r\in \B & (i,j)\in E, r\in \RR.
\end{align}
One of the imminent advantages of using objective $\zeta(\rho)$ is that we need only declare variables for entries $(i,j)\in E$ and can consequently delete the weak  constraints \eqref{MIP_lower_bool} from the formulation.
%
The LP relaxation of MIP($\rho$) (MLP($\rho$)) is obtained by giving up on the integrality constraints on $q_r$ and observing that without loss of generality we can simply  write $q_r\ge$ for all $r\in \RR$.
We next show that the optimal solutions of the LP relaxation of $\text{MIP}_{\text{F}}$ and MLP($\rho$) with $\rho=\frac{1}{k}$ coincide.
\begin{proposition}
	\label{lemma_MLPexact_equals_MLP1k}
	The optimal solutions of the LP relaxations $\textup{MLP}_{\textup{F}}$ and $\textup{MLP}(\frac{1}{k})$ coincide. 
\end{proposition}
\proof{Proof.} 
	It suffices to observe that as $\text{MLP}_{\text{F}}$ is a minimisation problem, each $z_{ij}$ $(i,j)\in \overline{E}$ takes the  value $\frac{1}{k}\sum_{ r\in \RR_{(i,j)} } q_r$  in any optimal solution to $\text{MLP}_{\text{F}}$ due to constraint \eqref{MIP_lower_bool}. 
	This implies that the second terms in the objective function \eqref{MIP_objective} of $\text{MIP}_{\text{F}}$ and \eqref{MIP_rho_objective} of MLP($\frac{1}{k}$) have the same value.\Halmos
\endproof
Therefore one may {instead solve MLP($\frac{1}{k}$) that has fewer variables and constraints than $\text{MLP}_{\text{F}}$.} 
In addition, for all $\rho>0$, a corollary of Proposition \ref{lemma_MLP_non_zero} holds by looking at the dual of MLP($\rho$) (MDP($\rho$)).
Let us associate variables $p_{ij}$ for $(i,j)\in E$ to constraints \eqref{MIP_rho_upper_bool} and  variable $\mu$ to constraint \eqref{MIP_rho_cardinality}. Then the dual of MLP($\rho$)  is:
\begin{align}
	(\text{MDP}(\rho))\;\;\zeta_{\text{MDP}(\rho)} =\max_{p,\mu} \;& \sum_{(i,j)\in E}p_{ij} - k \,\mu\\
	\text{s.t. }
	&\sum_{(i,j)\in E\cap \supp(R)  } p_{ij} -  \mu \le\rho\, |\overline{E} \cap \supp(R)|  & R\in \RR \label{MDP_rho_balance_con}\\
	& \mu\ge 0, \, p_{ij}\in[0,1]& (i,j)\in E
\end{align}
where $\supp (R) = \{(i,j) : r_{ij}=1  \}$.

\begin{corollary}
	Let $X$ have isolation number $i(X)>k$. Then for all $\rho>0$, \textup{MLP}($\rho$) has objective value at least $\rho \left(i(X) - k\right)$.
\end{corollary}
\proof{Proof.}
	The proof is a simple modification of Proposition \ref{lemma_MLP_non_zero}'s proof.
	The dual of MLP($\rho$) (MDP($\rho$)) differs from $\text{MDP}_{\text{F}}$ by having the constant value $\rho$ instead of dual variables $s_{ij}$ and constraints \eqref{MDP_rho_balance_con} instead of \eqref{MDP_balance_con}.
	Therefore setting $p_{ij}=\rho$ for all $(i,j)\in S$  and $0$ otherwise (where $S$ is a maximum isolated set of ones of $X$), and $\mu=\rho$ gives the required bound of $\rho \left(i(X) - k\right)$.\Halmos
\endproof

\section{Computational approach.}
\label{section_column_generation}

It is clearly not practical to solve the master integer program MIP($\rho$) or its LP relaxation MLP($\rho$) explicitly as the formulation has an exponential number of variables. 
\textit{Column generation} (CG) is a well-known technique to solve large LPs  iteratively by  only considering the variables which have the potential to improve the objective function \cite{Nemhauser:1998}. 
The column generation procedure is initialised by solving  a \textit{Restricted} Master LP (RMLP) which has a small subset of the variables of the full problem.
The next step is to identify a missing variable with \textit{negative reduced cost} to be added to RMLP.
To  avoid  considering all missing variables explicitly, a \textit{pricing problem} is formulated and solved. The solution of the pricing problem either returns a variable with negative reduced cost  and the procedure is iterated; or proves that no such variable exists and hence the solution of  RMLP is optimal for the full MLP.
In this section, we detail how CG technique can be used to solve the LP relaxation of MIP($\rho$) iteratively. 

Each \textit{Restricted} MLP($\rho$) (RMLP($\rho$)) has the same number of constraints as the full MLP($\rho$) and all variables $\xi_{ij}$ for $(i,j)\in E$ but it only has a small subset of variables $q_r$ for  $r\in\RR' \subset\RR$ where $|\RR'| \ll |\RR|$. 
Recall that each variable $q_r$ corresponds to a rank-$1$ binary matrix  $r\in\RR$ which determines the coefficients of $q_r$ in the constraints as well as the objective function. 
Hence at every iteration of the CG procedure we either need to find a rank-$1$ binary matrix for which the associated variable has a negative reduced cost, or, prove that no such matrix exists. 

\subsection{The pricing problem.}
At the first iteration of CG, RMLP($\rho$) may be initialised with $\RR'=\emptyset$ or can be warm started by identifying a few rank-$1$ matrices in $\RR$ using a heuristic.
After solving the RMLP($\rho$) to optimality, one  obtains an optimal dual solution  $[\mvec{p}^*,\mu^*]$ to the current RMLP($\rho$).
To identify a missing variable $q_r$ that has negative reduced cost, we solve the following pricing problem (PP): 
\begin{align}
(\text{PP})  \;\;\omega(\mu^*, \mvec{p}^*)&=\mu^* - \max_{a,b,y}  \sum_{(i,j)\in E} p_{ij}^{*} y_{ij} -  \rho \sum_{(i,j)\in \overline{E}} y_{ij} &\\
\text{s.t. }
& y_{ij}=a_ib_j, \qquad a_{i},b_{j}\in\B,\;&i\in [n], j\in [m].
\end{align}
PP may be formulated as an integer linear program $(\text{IP}_{\text{PP}})$ by using McCormick envelopes \cite{McCormick:1976} (see Section \ref{section_compact_formulation}) to linearise
the quadratic constrains  to $y_{ij}\in MC(a_i,b_j)$. 
The objective of PP depends on the current dual solution $[\mvec{p}^*,\mu^*]$ 
and its optimal solution corresponds to a rank-$1$ binary matrix $\mvec{a}\mvec{b}^\top=r\in\RR$ whose corresponding variable $q_r$ in MLP($\rho$) has the smallest reduced cost.
If $\omega( \mu^*,\mvec{p}^*)\ge 0$, then the current RMLP($\rho$)  does not have any missing variables with negative reduced cost and consequently  
the current solution of RMLP($\rho$) is optimal for MLP($\rho$).
If $\omega( \mu^*,\mvec{p}^*)<0$, then the variable $q_r$ associated with the rank-$1$ binary matrix $r=\mvec{a}\mvec{b}^\top$ is added to the next RMLP($\rho$) 
and the procedure is iterated.
Moreover, any feasible solution to PP with a negative reduced cost can (also) be added to the RMLP($\rho$) to continue the  procedure.
CG  terminates with a proof of optimality  if at some iteration we have $\omega(\mu^*, \mvec{p}^*)\ge 0$. 

\subsection{Solving the master integer program.}
After the CG process, if the optimal solution of MLP($\rho$) is integral, then it also is optimal for $\text{MIP}(\rho)$. However, if it is fractional, then this solution only provides a lower bound on the optimal value of $\text{MIP}(\rho)$.
In this case we obtain an integer feasible solution by solving a \textit{Restricted} MIP($\rho$) (RMIP($\rho$)) over the rank-$1$ binary matrices generated by the CG process applied to MLP($\rho$). 
This integer feasible solution is optimal for $\text{MIP}(\rho)$ provided that the objective value of  RMIP($\rho$) is equal to the ceiling of the objective value of MLP($\rho$).
If this is not the case, one needs to embed CG into a branch-and-bound tree \cite{Lubbecke:2005} to solve MIP($\rho$) to optimality,  which is a relatively complicated process and we do not consider it in this paper. 

\subsection{Computing lower bounds.}
Note that even if the CG procedure is terminated prematurely, one can still obtain a lower bound on MLP($\rho$) and therefore on $\text{MIP}(\rho)$ by considering the dual of MLP($\rho$).
Let the objective value of  of the current RMLP($\rho$) be
\begin{align}
\zeta_{\text{RMLP}(\rho)}= \sum_{(i,j)\in E} \xi^*_{ij} +\rho \sum_{(i,j) \in \overline{E}} \sum_{r\in \RR_{(i,j)}} q^*_r
= \sum_{(i,j)\in E}p^*_{ij} - k \mu^*
\end{align}
where  $[\xi^*_{ij}, q^*_r]$  is the optimal solution of  RMLP($\rho$) and  $[\mvec{p}^*, \mu^*]$ is the corresponding optimal dual solution which does not necessarily satisfy all of the constraints \eqref{MDP_rho_balance_con} for MDP($\rho$).
Now assume that we solve PP to optimality and  obtain a rank-$1$ binary matrix with a negative reduced cost, $\omega(\mu^*, \mvec{p}^*)<0  $. In this case, we can construct a feasible solution $[\mvec{p}, \mu]$ to MDP($\rho$) by setting $\mvec{p}:=\mvec{p}^*$ and $\mu :=\mu^*- \omega(\mu^*, \mvec{p}^*)$ and obtain the following bound on the optimal value $\zeta_{\text{MLP}(\rho)}$ of MLP($\rho$),
\begin{equation}
\label{guarantee}
\zeta_{\text{MLP}(\rho)} \ge\sum_{(i,j)\in E}p_{ij}  - k\,\mu
= \sum_{(i,j)\in E}p^*_{ij} - k\,\left( \mu^*- \omega( \mu^*,\mvec{p}^*) \right) 
= \zeta_{\text{RMLP}(\rho)} + k \; \omega( \mu^*,\mvec{p}^*).
\end{equation}
If we do not have the optimal solution to PP but have a lower bound $\underline \omega(\mu^*, \mvec{p}^*)$ on it, $\omega( \mu^*,\mvec{p}^*)$ can be replaced by $\underline \omega(\mu^*, \mvec{p}^*)$ in Equation \eqref{guarantee} and the bound on MLP($\rho$) still holds. Furthermore, this lower bound on MLP($\rho$) naturally provides a valid lower bound on $\text{MIP}(\rho)$, thus giving us a bound on the optimality gap.

\subsection{Column generation for \texorpdfstring{$\text{MLP}_{\text{F}}$}{MLP F}. }
The CG approach  is described above as applied to the LP relaxation of $\text{MIP}(\rho)$. 
To apply CG to $\text{MLP}_{\text{F}}$ only a small modification needs to be done.
The Restricted $\text{MLP}_{\text{F}}$ provides dual variables for constraints \eqref{MIP_lower_bool} which are used in the objective of PP for coefficients of $y_{ij}$ $(i,j) \in \overline{E}$. 

We note that CG cannot be used to solve the LP relaxation of the \textit{strong formulation} of $\text{MIP}_{\text{F}}$ in which constraints  \eqref{MIP_lower_bool} are replaced by exponentially many constraints $q_r\le z_{ij}$ for all $r\in \RR_{(i,j)}$ and $(i,j)\in \overline{E}$.
This is due to the fact that CG could cycle and generate the same column over and over again.
For example, consider applying CG to solve the strong formulation of $\text{MLP}_{\text{F}}$ and start with the rank-$1$ binary matrix of all $1$s as the first column associated with variable $q_1$. 
The objective value of the corresponding Restricted $\text{MLP}_{\text{F}}$  would be $\zeta_{\text{RMLP}}^{(1)} = 0 + |\overline{E}|$ for the solution vector $[\mvec{\xi}^{(1)}, \mvec{z}^{(1)},\mvec{q}^{(1)}] = [ \mvec{0}, \mvec{1}, 1]$ as all entries of the input matrix are covered. 
Adding the same rank-$1$ binary matrix of all $1$s in the next iteration and setting $[q_1,q_2] = [\frac{1}{2}, \frac{1}{2}]$, allows us to keep $\mvec{\xi}^{(2)} = \mvec{0}$ but reduce the value of  $\mvec{z}^{(2)} $ to $ \frac{1}{2}\mvec{1}$ to obtain an objective value $\zeta_{\text{RMLP}}^{(2)} = 0 + \frac{1}{2}|\overline{E}|$. 
Therefore, repeatedly adding the same matrix of all $1$s for $t$ iterations, the objective function would become $\zeta_{\text{RMLP}}^{(t)} = 0 + \frac{1}{t}|\overline{E}|$ for the solution vector $[\mvec{\xi}^{(t)}, \mvec{z}^{(t)},\mvec{q}^{(t)}] = [\mvec{0}, \frac{1}{t}\mvec{1}, \frac{1}{t}\mvec{1}]$. 
Consequently,  as $t\to \infty$ we would have $\zeta_{\text{RMLP}}^{(t)} \to 0$ and during the column generation process we repeatedly generate the same rank-1 binary matrix.

\subsection{An alternative formulation of the pricing problem.}
Generating rank-$1$ binary matrices with negative reduced cost efficiently is at the heart of the CG process. 
For both MLP($\rho$) and $\text{MLP}_{\text{F}}$,  the pricing problem is a Bipartite Binary Quadratic Program (BBQP) which is NP-hard in general \cite{Punnen:2012,Peeters:2003}. Hence for large $X$ it may take too long to solve PP to optimality via formulation $\text{IP}_{\text{PP}}$ at each iteration.
Introducing  $H$ an $n\times m$ matrix with $h_{ij} = p^*_{ij}\in [0,1]$ for $(i,j)\in E$, $h_{ij} = -\rho$ for $(i,j) \in \overline{E}$ and $h_{ij}=0$ for $(i,j)\not \in E \cup \overline{E}$, PP can be written in standard form as 
\begin{equation}
\label{equation_PP_BBQ}
(\text{QP}_{\text{PP}}) \quad \omega(\mu^*,\mvec{p}^*) =\mu^*- \max_{\mvec{a}\in \B^n, \mvec{b} \in \B^m} \mvec{a} ^\top H \mvec{b}.
\end{equation}
This explicit quadratic form $\text{QP}_{\text{PP}}$ is more intuitive for thinking about heuristics than formulation $\text{IP}_{\text{PP}}$. 
If a heuristic approach to PP returns a rank-$1$ binary matrix with negative reduced cost, 
then it is valid to add this heuristic solution as a column to the next RMLP. 
\cite{Punnen:2012} presents several heuristics for BBQP along with a simple but powerful greedy algorithm. In Appendix \ref{appendix_pricing} we detail this greedy algorithm and some variants of it which we use to provide a warm start to PP at every iteration of CG in Section \ref{subsection_cg_test}.

\section{Experiments.}
\label{section_experiments}

The integer programs and column generation approach introduced in the previous sections provide a framework for computing $k$-BMF with dual bounds. In this section, we present some experimental results to demonstrate the practical applicability of integer programming to obtain low-error factorisations. More specifically we detail our pricing strategies during the column generation process and present a thorough comparison of models $\text{MIP}_{\text{F}}$, MIP($\rho$) and CIP on synthetic and real world datasets. Our code and data can be downloaded from \cite{Kovacs:code:2021}.

\subsection{Data.}

If $X$ contains rows (or columns) of all zeros, deleting these rows (or columns) leads to an equivalent problem whose solution $A$ and $B$ can easily be translated to a solution for the original problem by inserting a row of zeros to $A$ (respectively a column of zeros to $B$) in the corresponding place.
In addition, if $X$ contains duplicate rows or columns, by Lemma \ref{lemma_preprocessing} there is an optimal rank-$k$ factorisation which has the same row-column repetition pattern as $X$. 
Hence we solve the problem on a smaller matrix $X'$ which is obtained from $X$ by keeping  only one copy of each row and column, and use an updated objective function in which every entry is weighted proportional to the number of rows and columns it is contained in $X$.

\subsubsection{Synthetic data.}
\label{section_artificial_data}

We build our dataset of binary matrices with prescribed sparsity and Boolean rank as follows. To get a matrix $X\in \B^{n\times m}$ with Boolean rank at most $\kappa$, first we randomly generate two binary matrices $\tilde{A}$, $\tilde{B}$ of dimension $n\times \kappa$ and $\kappa\times m$, then compute their Boolean product to get $X$. This ensures $X$ has Boolean rank at most $\kappa$.
To obtain a certain sparsity for $X$, we control the probability of entries of $\tilde{A}$, $\tilde{B}$ being zero.
More specifically, if we generate $\tilde{a}_{i\ell}$, $\tilde{b}_{\ell j}$ to be zero with probability $p$, then $x_{ij}=\bigvee_{\ell=1}^{\kappa} \tilde{a}_{i\ell} \tilde{b}_{\ell j}$ is zero with probability $(1-(1-p)^2)^{\kappa}$.
Hence, to obtain $X$ with $\sigma$ percent of zeros, 
we need to generate entries of $\tilde{A}$, $\tilde{B}$ to be zero with probability 
$p=1-\sqrt{1-(\sigma/100)^{\frac{1}{\kappa}}}$. 

We generate matrices as described above with $m=20$ columns and  $\kappa=10$. The number of rows ($n$) is set to be $20, 35$ or $ 50$. For each of the three dimensions ($20\times20, 35\times 20, 50 \times 20$), we generate 10 \textit{sparse} matrices with $75$\% zeroes and 10 \textit{normal} matrices with $50$\% zeroes, corresponding to 10 different seed settings in the random number generation.
We call this initial set of $2\cdot 3 \cdot 10$ matrices the \textit{clean} matrices.
Next, we create a set of \textit{noisy} matrices from the clean matrices  by randomly flipping $5\%$ of the entries of each matrix. The noisy matrices are not necessarily of Boolean rank at most $\kappa=10$, but they are at most $0.05 \cdot n\cdot m$ squared Frobenius distance away from a Boolean rank $10$ matrix. Therefore, our test bed consists of $120$ matrices corresponding to $2$ noise level settings (\textit{noisy} or \textit{clean}), $2$ sparsity levels (\textit{sparse} or \textit{normal}), $3$ dimensions ($20\times20, 35\times 20, 50 \times 20$) and $10$ random seeds. 
Applying the preprocessing steps to our synthetic dataset achieves the largest dimension  reduction on clean matrices, while the dimension of noisy matrices scarcely changes.
A table summarising the parameters used to generate our data can be found in the Appendix \ref{appendix_synthetic_data}. 

\subsubsection{Real world data.}
We work with eight real world categorical datasets that were downloaded from online repositories \cite{UCI, BOOKS}.
In general if a dataset has a categorical feature $C$ with $N$ discrete options $v_j$, $(j\in [N])$, we convert feature $C$ into $N$ binary features $B_j$ $(j\in N)$ so that if the $i$-th sample takes option $v_j$ for $C$ that is $(C)_i = v_j$, then we have  $(B_j)_i = 1$ and $(B_\ell)_i = 0$ for all $\ell\not = j \in [N]$. This technique of binarisation of categorical columns has been applied in \cite{Kovacs:2017} and \cite{Barahona:2019}. 
If a row $i$ has a missing value in the column of feature $C$, we leave the corresponding binary feature columns with missing values in row $i$.
Table \ref{dimensions} shows a short summary of the resulting full-binary datasets used, in-depth details on converting categorical columns into binary, missing value treatment and feature descriptions can be found in Appendix \ref{appendix_real_data}.
\begin{table}[ht]
	\centering
	\caption{\label{dimensions}Summary of binary real world datasets}
		\begin{tabular}{lcccccccc}
			&zoo & tumor 
			& hepatitis
			&heart& lymp & audio& apb& votes\\
		\hline
\up\down
			$n\times m$    & 101 $\times$ 17   & 339 $\times$ 24 & 155$\times$ 38   & 242$\times$ 22  & 148$\times$44   & 226 $\times$ 92  & 105$\times$ 105  & 435 $\times$16\\
			\# missing & 0 & 670 & 334 & 0&0 & 899& 0&392 
			\\
			\%1s  & 44.3  & 24.3 & 47.2  & 34.4 & 29.0  & 11.3  & 8.0  & 49.2 \\
		\hline
		\end{tabular} %
\end{table}

\subsection{Testing the computational approach to exponential formulation I.}
\label{subsection_cg_test}

Since the efficiency of CG greatly depends on the speed of generating columns, let us illustrate the speed-up gained by using heuristics to solve the pricing problem. 
At each iteration of CG, 30 heuristic solutions are computed via the heuristics detailed in Appendix \ref{appendix_pricing} in order to obtain initial feasible solutions to PP.
%
Under \textit{exact} pricing, the best heuristic solution is used as a warm start and $\text{IP}_{\text{PP}}$ is solved to optimality at each iteration using CPLEX \cite{CPLEXmanual}.
In simple heuristic (\textit{heur}) pricing, if the best heuristic solution to PP has negative reduced cost 
then it is directly added to the next RMLP($\rho$). If at some iteration, the best heuristic column does not have negative reduced cost, CPLEX is used to solve $\text{IP}_{\text{PP}}$ to optimality for that iteration.
The multiple heuristic (\textit{heur\_multi}) pricing strategy is a slight modification of the simple heuristic strategy, in which at each iteration all columns with negative reduced cost are added to the next RMLP($\rho$). 

\begin{figure}[h]
		\caption{\label{zoo_pricing}Comparison of pricing strategies for solving $\text{MLP}(1)$ on the zoo dataset} 
	\centering
	\begin{subfigure}[b]{0.4\textwidth}
		\centering
		\includegraphics[scale=0.65]{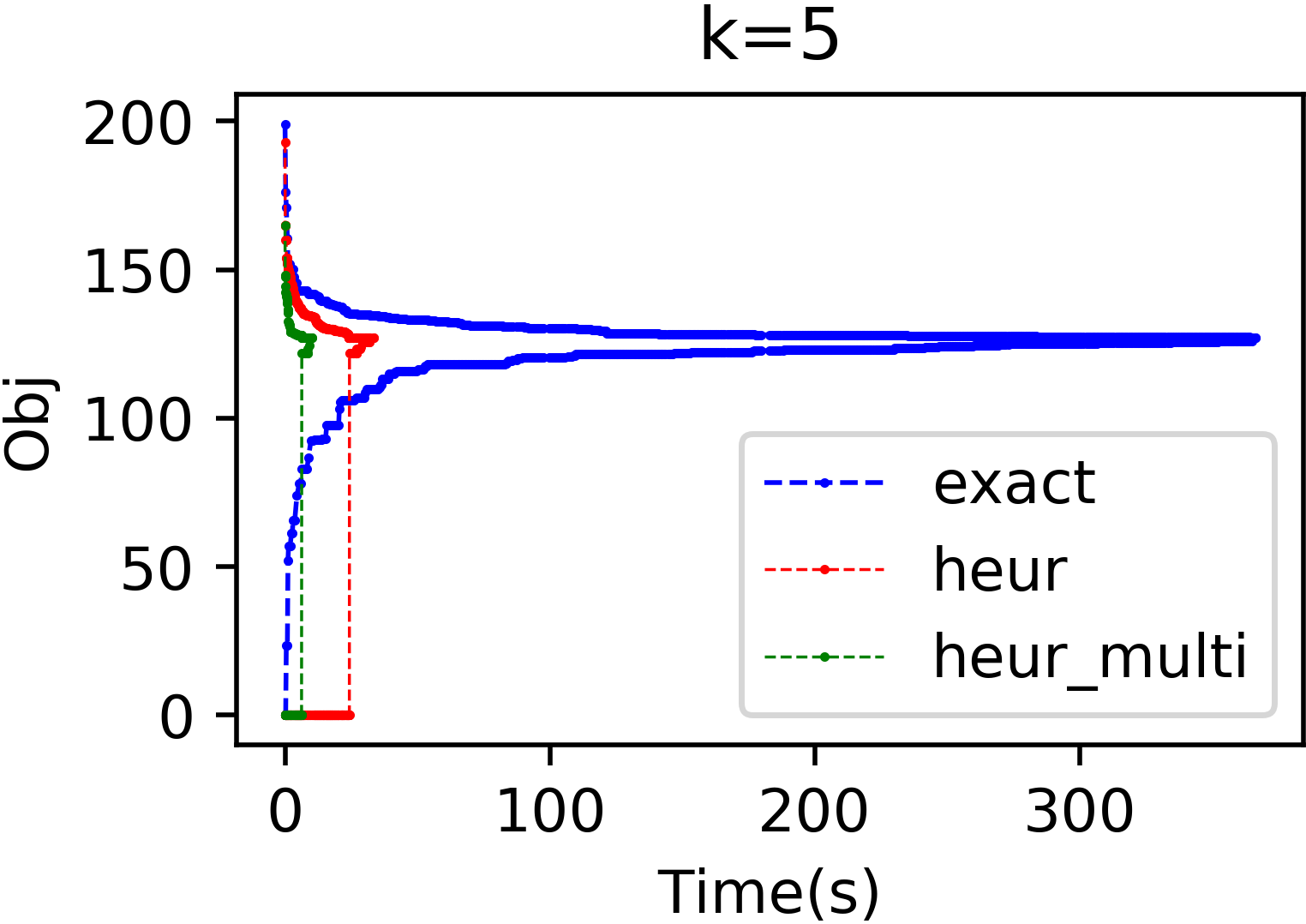}
	\end{subfigure}
	\begin{subfigure}[b]{0.4\textwidth}
		\centering
		\includegraphics[scale=0.65]{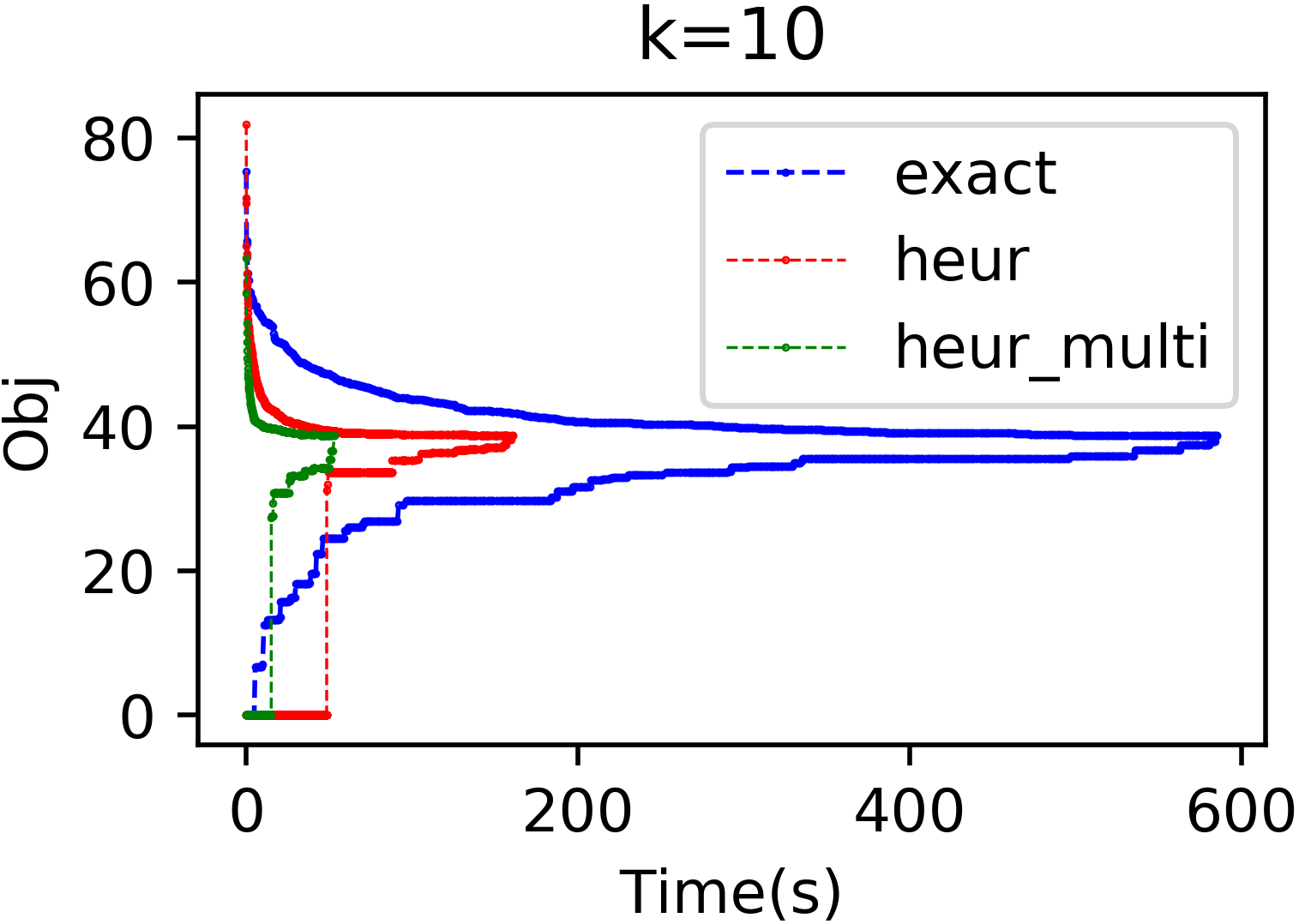}
	\end{subfigure}
\end{figure}

Figure \ref{zoo_pricing} indicates the differences between pricing strategies when solving MLP(1) via CG for $k=5,10$ on the zoo dataset. The primal objective value of MLP(1) (decreasing curve) and the value of the dual bound (increasing curve) computed using the formula in Equation \eqref{guarantee} are plotted against time. 
Sharp increases in the dual bound for heuristic pricing strategies correspond to iterations in which CPLEX was used to solve $\text{IP}_{\text{PP}}$, as for the evaluation of the dual bound on MLP(1) a lower bound on ${\omega}(\mu^*,\mvec{p}^*)$ is needed which heuristic solutions do not provide. 
While we observe a tailing off effect \cite{Lubbecke:2005} on all three curves, both heuristic pricing strategies provide a significant speed-up from exact pricing, adding multiple columns at each iteration being the fastest.

In order for CG to terminate with a certificate of optimality, at least one pricing problem has to be solved to optimality. Unfortunately for larger datasets we cannot expect this to be achieved in a short amount of time. 
Therefore, we change the multiple heuristic pricing strategy to get a pricing strategy that we use in the rest of the experiments as follows.
We impose an overall fixed time limit  on the CG process and use the barrier method in CPLEX as the LP solver for RMLP at each iteration. 
At each iteration of CG, we add up to 2 columns with the most negative reduced cost to the next RMLP.
If at an iteration, heuristics for PP do not provide a column with negative reduced cost and CPLEX has to be used to improve the heuristic solution, we do not solve $\text{IP}_{\text{PP}}$  to optimality but abort CPLEX after 25 seconds if a column with negative reduced cost has been found.
While these modifications result in a speed-up, they reduce the chance of obtaining a strong dual bound. 
In case we wish to focus more on computing a stronger dual bound on MLP,  
we may continue solving $\text{IP}_{\text{PP}}$ via CPLEX even when a heuristic negative reduced cost solution is available.

\subsubsection{MLP(1) vs \texorpdfstring{$\text{MLP}_{\text{F}}$}{MLP F}. }

In this section we compare the LP relaxations of MIP(1) and $\text{MIP}_{\text{F}}$. 
According to Proposition \ref{lemma_MLPexact_equals_MLP1k} the optimal solution of  $\text{MLP}_{\text{F}}$ is equivalent to MLP($\frac{1}{k}$) and hence we solve MLP($\frac{1}{k}$) which has fewer variables and constraints than $\text{MLP}_{\text{F}}$. 
To solve MLP(1) and MLP($\frac{1}{k}$), we start off from $0$ rank-$1$ binary matrices so $\RR'=\emptyset$ in the first RMLP and  
set a total time limit of $600$ seconds, so we either solve MLP to optimality under 600 seconds or run out of time and compute the gap between the last RMLP and the best dual bound MDP according to formula $100(\zeta_{\text{RMLP}}-\zeta_{\text{MDP}})/\zeta_{\text{RMLP}} $. 
As MLP(1) and MLP($\frac{1}{k}$) correspond to the LP relaxations of MIP(1) and  $\text{MIP}_{\text{F}}$ with integral objective coefficients, any fractional dual bound may be rounded up to give a valid bound on the master IP.
Therefore, we stop CG whenever the ceiling of the dual bound 
 reaches the objective value of RMLP.
\begin{figure}[htbp]
		\caption{\label{figure_MLP1_MLP1k_time_cols} Time taken in seconds to solve MLP(1) and MLP($\frac{1}{k}$) via CG on synthetic data }
	\centering
	\begin{subfigure}[b]{0.24\textwidth}
		\includegraphics[scale=0.48]{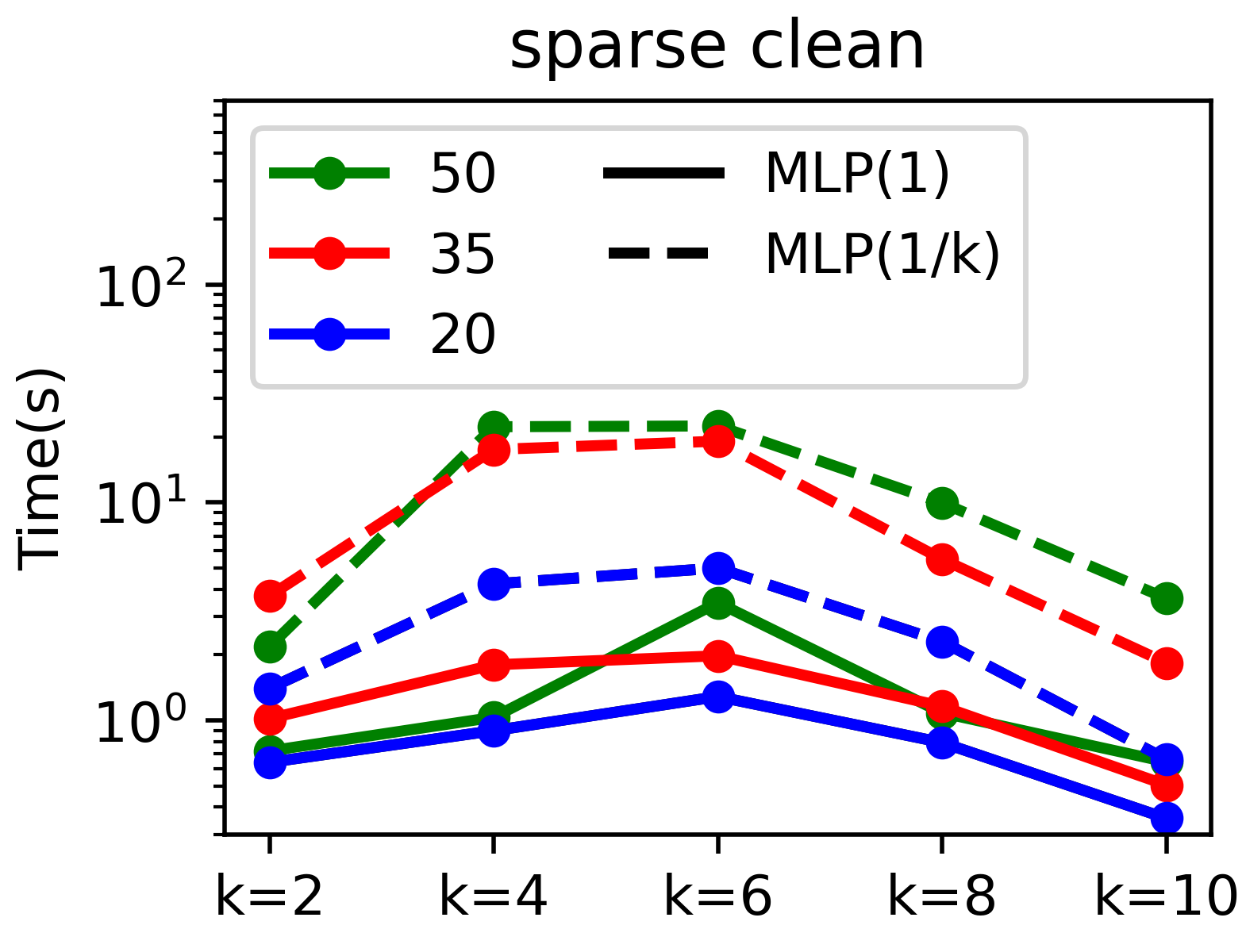}
		\label{subfigure_MLP1_MLP1k_time_sparse_clean}
	\end{subfigure}
	\hfill
	\begin{subfigure}[b]{0.24\textwidth}
		\includegraphics[scale=0.48]{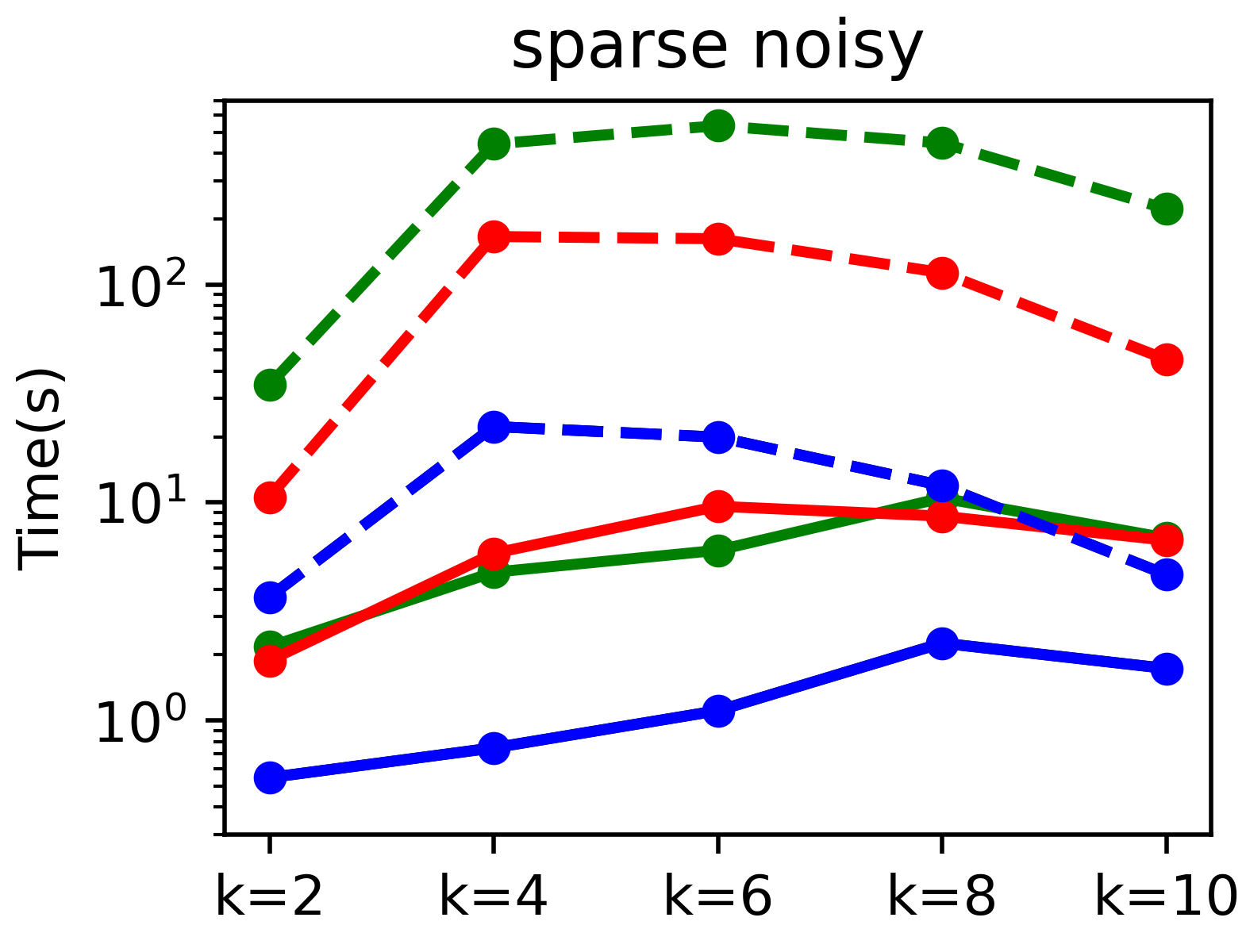}
		\label{subfigure_MLP1_MLP1k_time_sparse_noisy}
	\end{subfigure}
	\hfill
	\begin{subfigure}[b]{0.24\textwidth}
		\includegraphics[scale=0.48]{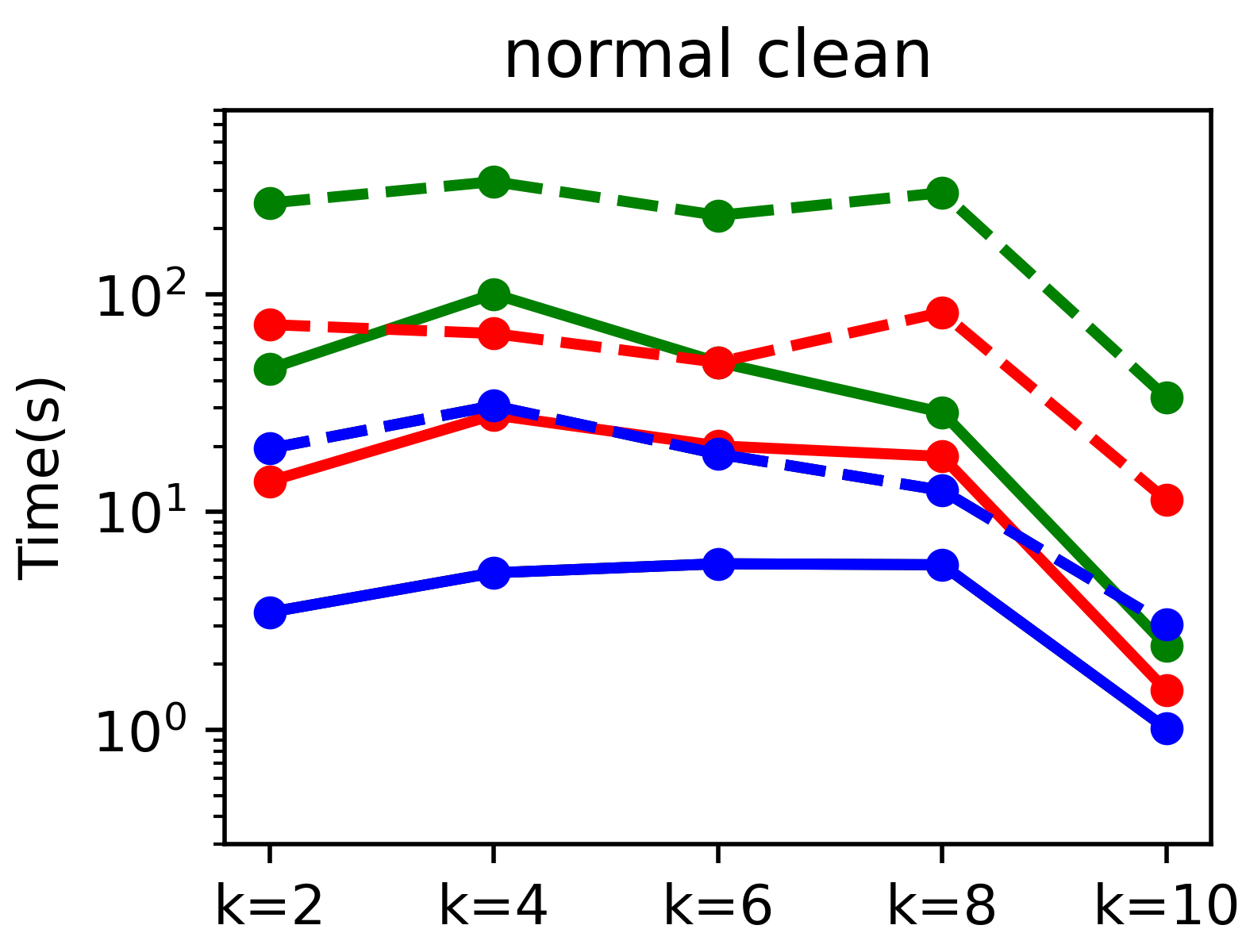}
		\label{subfigure_MLP1_MLP1k_time_normal_clean}
	\end{subfigure}	
	\hfill	
	\begin{subfigure}[b]{0.24\textwidth}
		\includegraphics[scale=0.48]{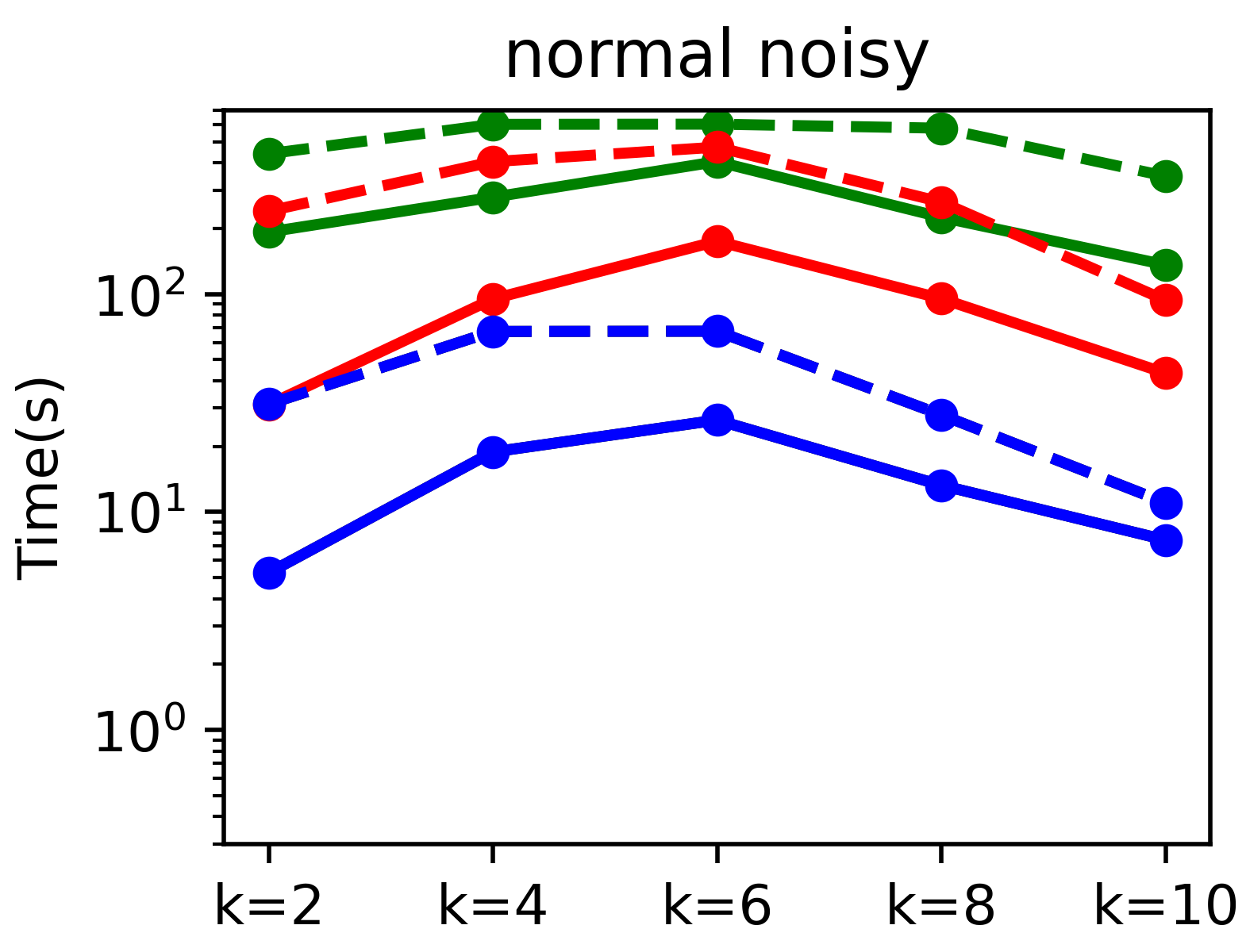}
		\label{subfigure_MLP1_MLP1k_time_normal_noisy}
	\end{subfigure}
\end{figure}

Figure \ref{figure_MLP1_MLP1k_time_cols} shows the time taken in seconds on a logarithmic scale to solve MLP(1) and MLP($\frac{1}{k}$) via CG for $k=2,4,\dots,10$ on the synthetic matrices. 
Each line corresponds to the average  taken over $10$ instances with the same dimension, sparsity and noise level. Blue lines correspond to matrices of dimension $20\times 20$, red to $35\times 20$ and green to $50 \times 20$. Solid lines are used for MLP(1) and dashed for MLP($\frac{1}{k}$).
%
First, we observe that it is significantly faster to solve both MLPs on \textit{sparse} and \textit{clean} matrices as opposed to \textit{normal} and \textit{noisy} ones of the same dimension.
Preprocessing is more effective in reducing the dimension for clean matrices in comparison to noisy ones (see Table \ref{table_synthetic_data} in Appendix \ref{appendix_synthetic_data}) which explains why noisy instances take longer. In addition, both MLP(1) and MLP($\frac{1}{k}$) have a number of variables and constraints directly proportional to non-zero entries of the input matrix, hence a sparse input matrix requires a smaller problem to be solved.
%
Second, we see that $k=10$ are solved somewhat faster. 
This can be explained by all matrices in our test bed being generated to have Boolean rank at most $10$. 
For a rank-$10$ factorisation of \textit{clean} matrices without noise we get $0$ factorisation error under both models MIP(1) and $\text{MIP}_{F}$ and hence LP relaxation objective value $0$.
For \textit{noisy} matrices we observe the error to be in line with our expectation of  $0.05 \cdot n \cdot m$.
%
We observe that in some cases it takes significantly longer to solve MLP($\frac{1}{k}$), and in all ten instances of $50\times 20$ \textit{normal}-\textit{noisy} matrices MLP($\frac{1}{k}$) for $k=6$ runs out of the time budget of $600$ sec. In the experiments, we see the amount of time CG takes is directly proportional to the number of columns generated, MLP($\frac{1}{k}$) generating significantly more columns than MLP(1).

\subsubsection{Obtaining integral solutions.}
\label{subsection_obtaining_integer_sols}
Once we obtain some rank-$1$ binary matrices (i.e. columns) via CG applied to a master LP, we can obtain an integer feasible solution by solving either of the master IPs over the columns available. Here we explore obtaining integer feasible solutions by solving MIP(1) and $\text{MIP}_{F}$ over the columns generated by formulations MLP(1) and MLP($\frac{1}{k}$). We use CPLEX as our integer program solver and set a total time limit of $300$ seconds.
\begin{figure}[htbp]
		\caption{\label{figure_MLP1_MLP1k_MIP1_obj} Factorisation error in $\|\cdot \|_F^2$ of integral solutions by MIP(1) from columns  by MLP(1) and MLP($\frac{1}{k}$) }
	\centering
	\begin{subfigure}[b]{0.24\textwidth}
		\includegraphics[scale=0.48]{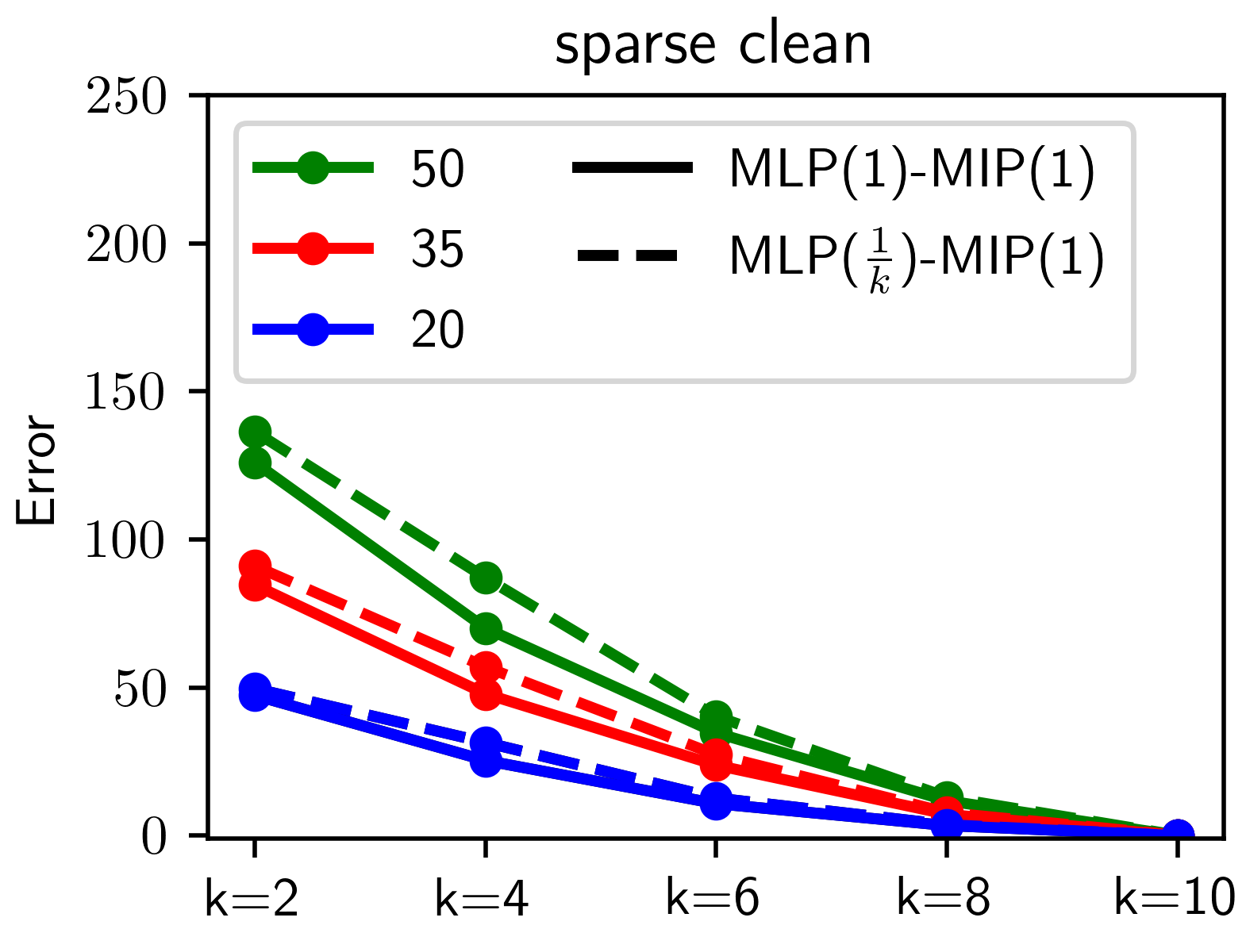}
	\end{subfigure}
	\hfill
	\begin{subfigure}[b]{0.24\textwidth}
		\includegraphics[scale=0.48]{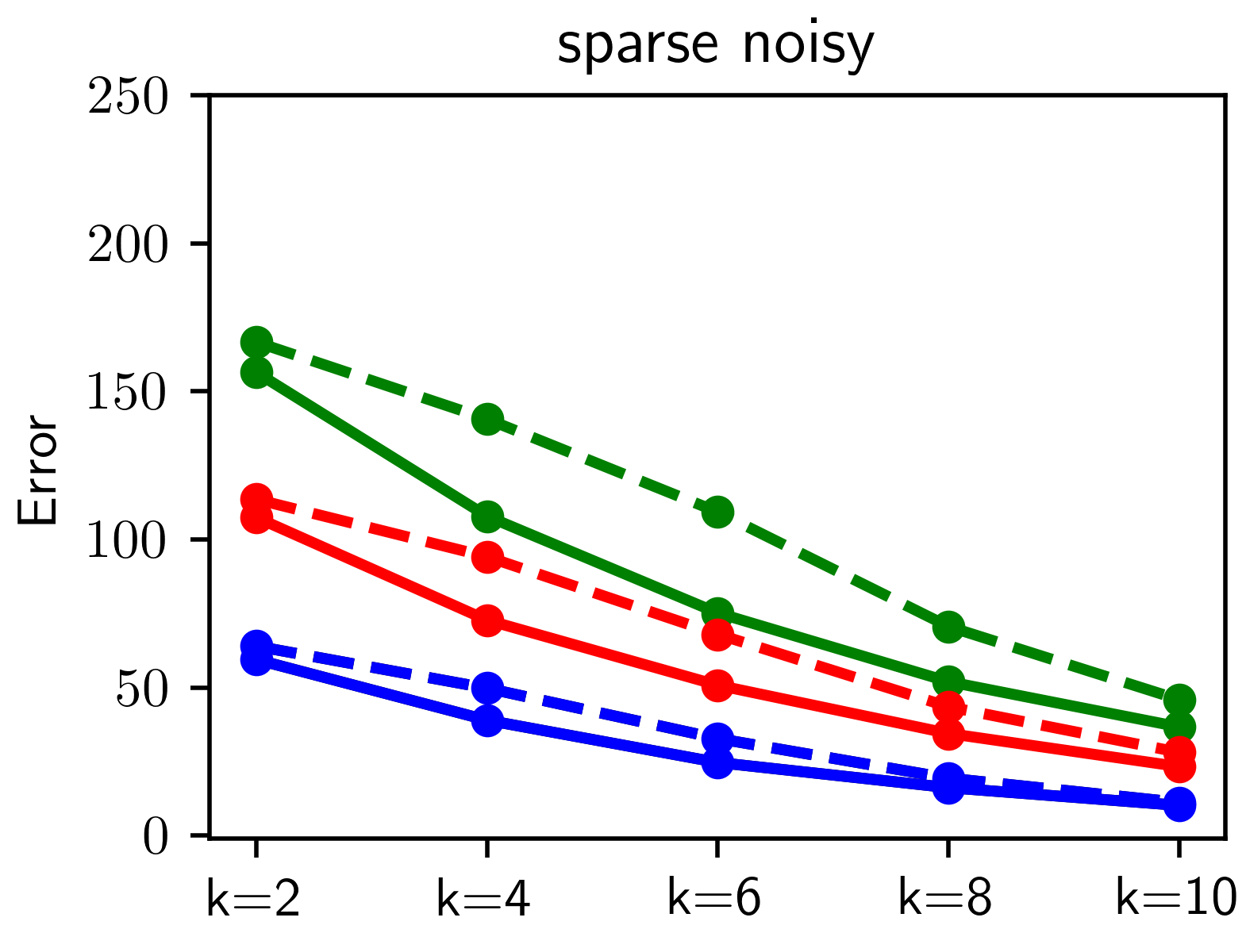}
	\end{subfigure}
	\hfill
	\begin{subfigure}[b]{0.24\textwidth}
		\includegraphics[scale=0.48]{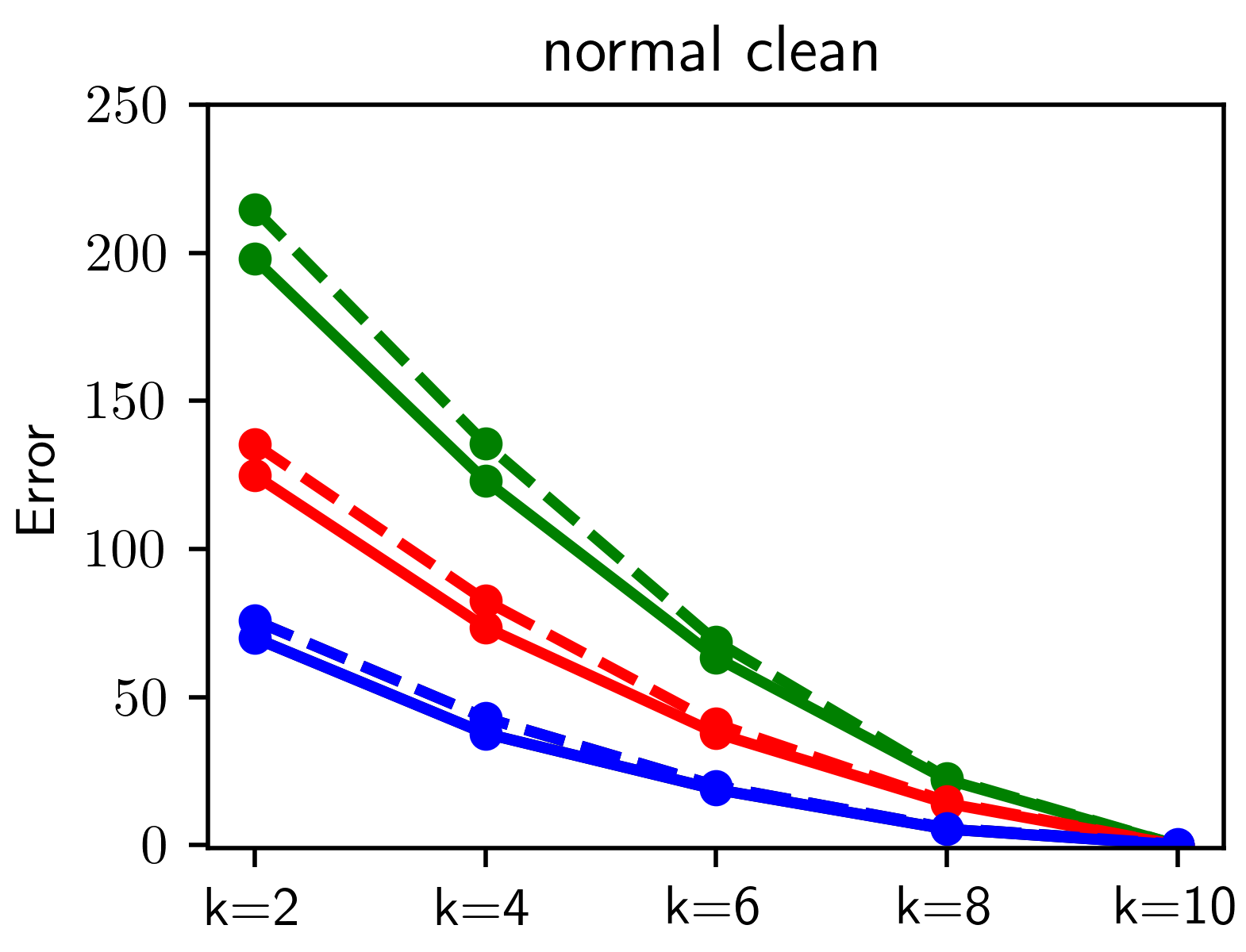}
	\end{subfigure}	
	\hfill	
	\begin{subfigure}[b]{0.24\textwidth}
		\includegraphics[scale=0.48]{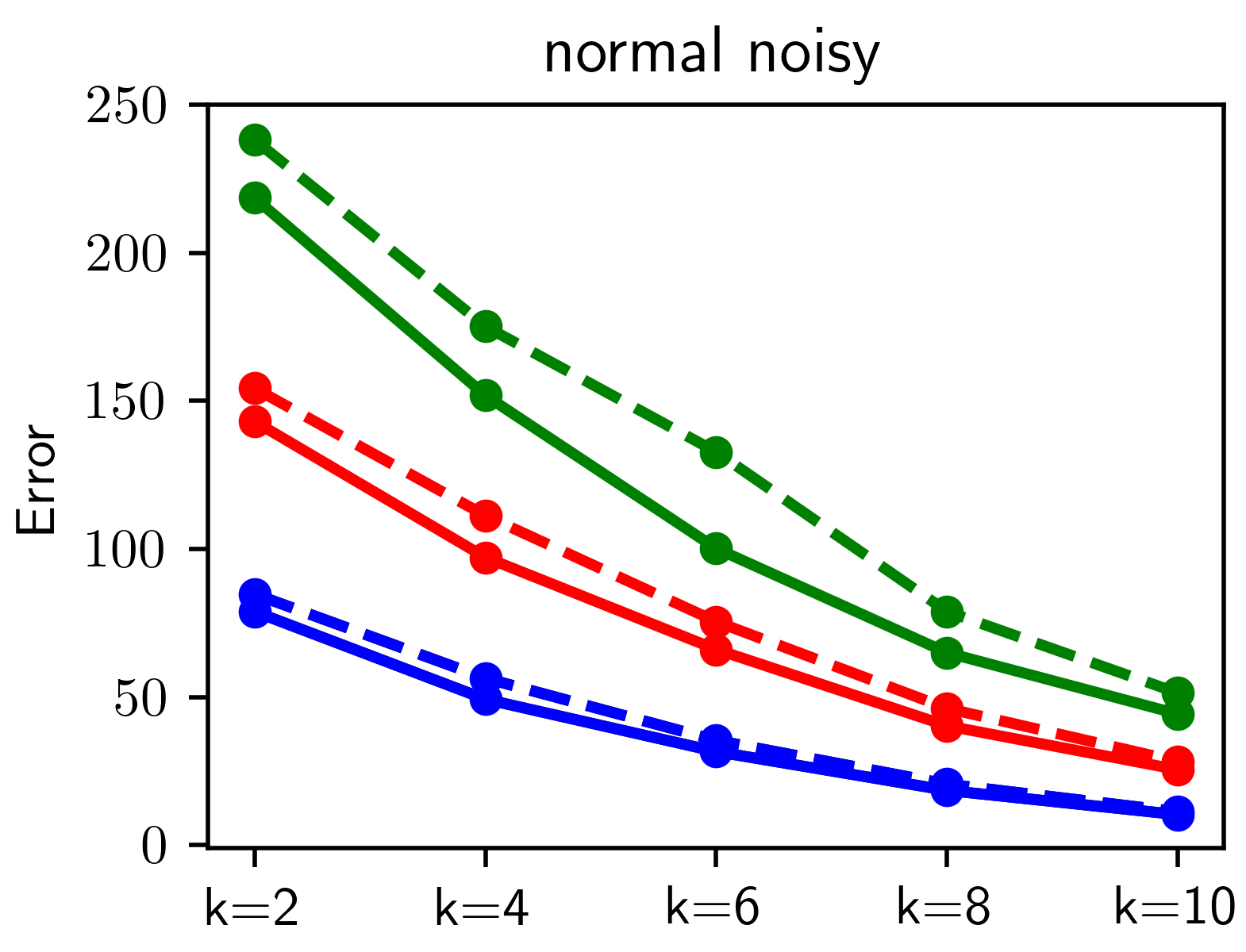}
	\end{subfigure}
\end{figure}

Figure \ref{figure_MLP1_MLP1k_MIP1_obj} shows the factorisation error  in $\|\cdot \|_F^2$ of integer feasible solutions obtained by solving MIP(1) over columns generated by MLP(1) and MLP($\frac{1}{k}$). As previously, each line corresponds to the average taken over 10 matrices with same dimension, sparsity and noise level. 
Solid lines are used to denote where the columns used were generated by MLP(1) and dashed where by MLP($\frac{1}{k}$).
Comparing the error values of the dashed and solid lines we draw a crucial observation: columns generated by MLP(1) seem to be a better basis for obtaining low-error integer feasible solutions than columns by MLP($\frac{1}{k}$).
We suspect this is the case as in the majority of rank-$k$ factorisations most entries are only covered by a few rank-$1$ binary matrices whereas MLP($\frac{1}{k}$) favours rank-$1$ matrices which heavily cover $0$ entries of the input matrix. 
This is because the coefficient in MLP($\frac{1}{k}$)'s objective function corresponding to a zero entry at position $(i,j)$  is only $\frac{1}{k} \times$(number of rank-$1$ matrices covering $(i,j)$), hence it is cheaper  for MLP($\frac{1}{k}$) to cover a $0$ by a few (less than $k$) rank-$1$ matrices than to leave any $1$s uncovered.
We also conducted a set of experiments using formulation $\text{MIP}_{\text{F}}$ and we see that the factorisation error when using formulation MIP(1) to obtain the integral solutions is extremely close to that of $\text{MIP}_{\text{F}}$, see Appendix \ref{appendix_obtaining_integer_sols} Tables \ref{table_MIP1_MIPexact_from_MLP1} and \ref{table_MIP1_MIPexact_from_MLP1k} for the precise difference in the factorisation error between the two master IPs.
\begin{figure}[htbp]
		\caption{\label{figure_MLP1_MIPs_time} Time taken in seconds to solve MIP(1) and $\text{MIP}_{F}$ on columns generated by MLP(1)}
	\centering
	\begin{subfigure}[b]{0.24\textwidth}
		\includegraphics[scale=0.48]{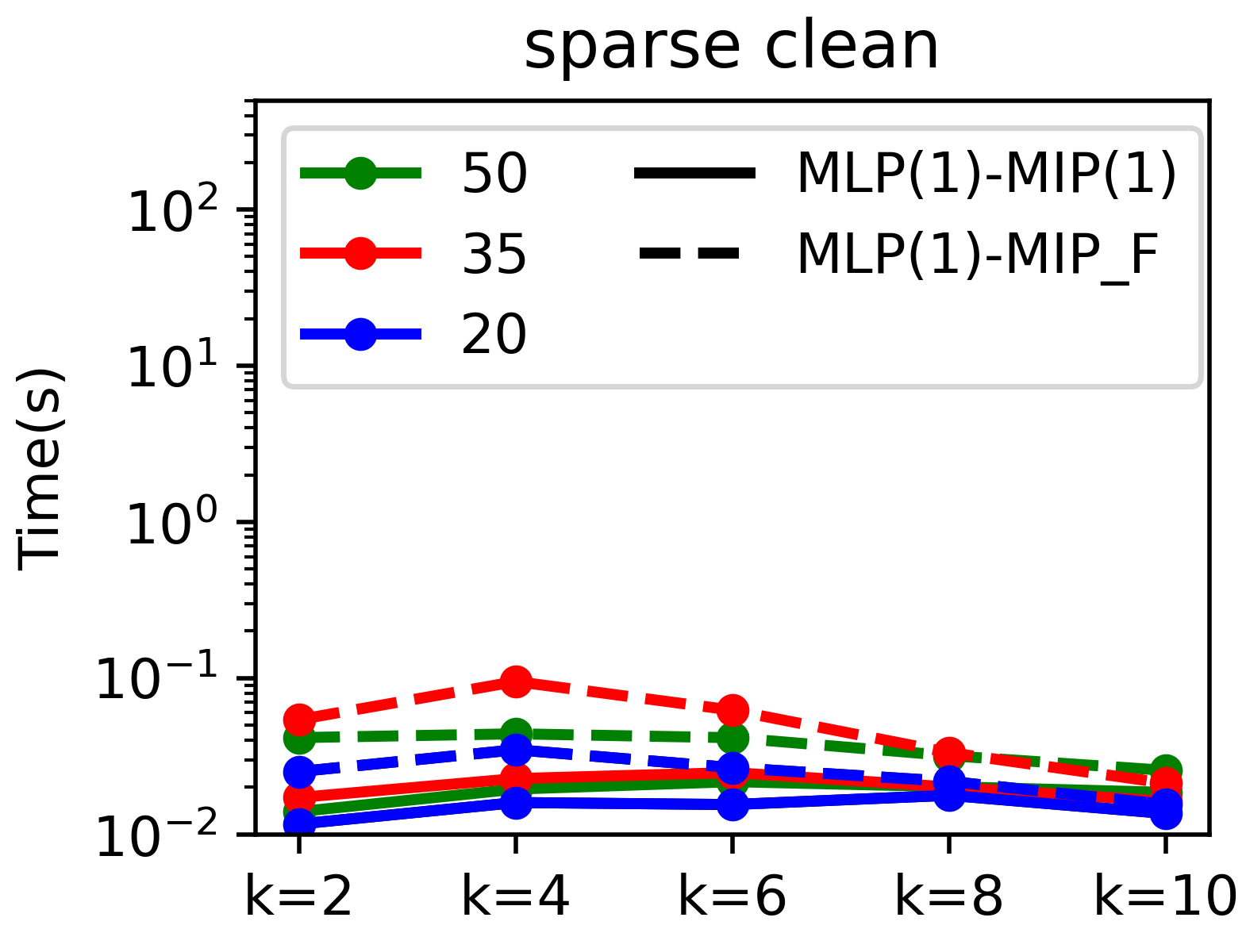}
	\end{subfigure}
	\hfill
	\begin{subfigure}[b]{0.24\textwidth}
		\includegraphics[scale=0.48]{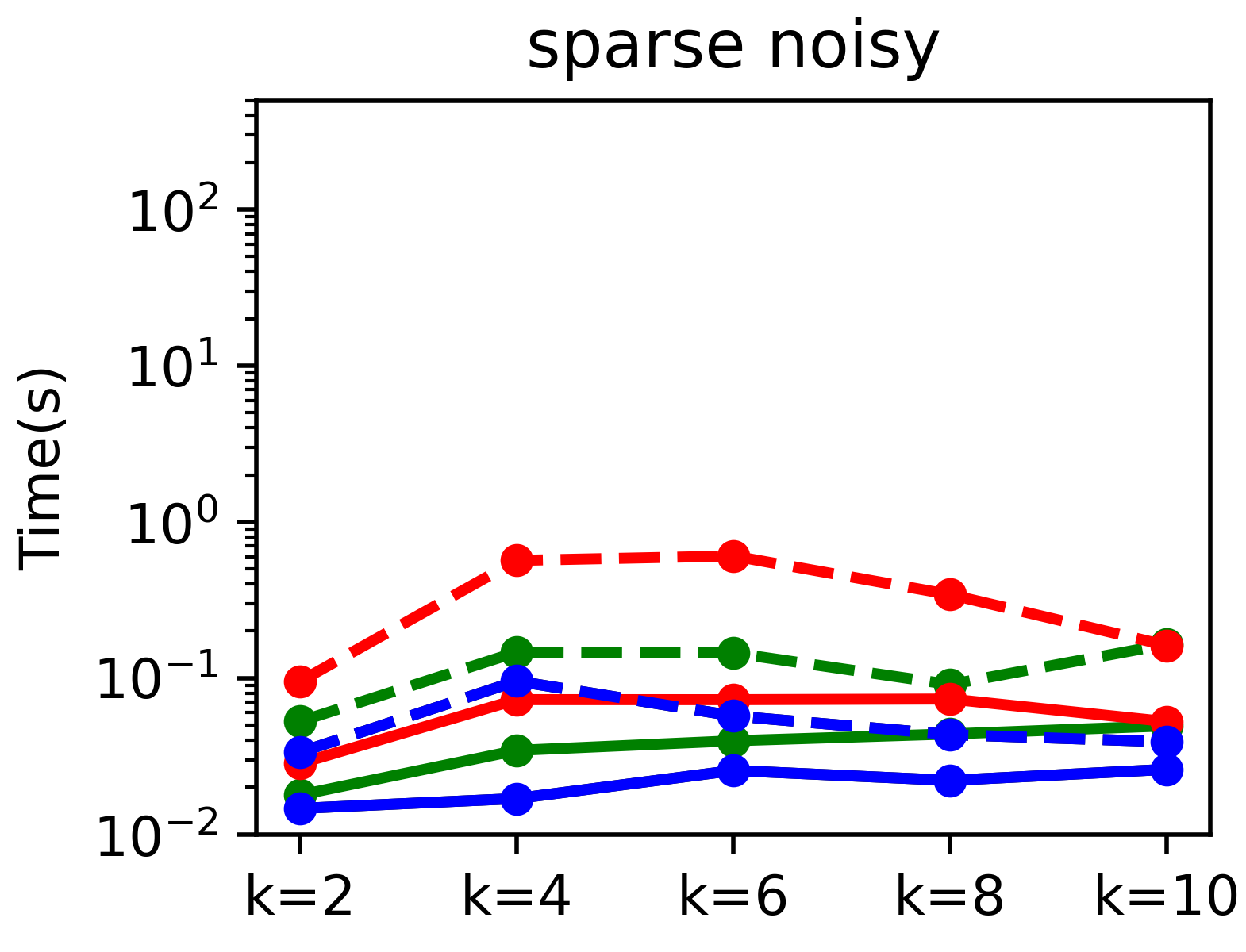}
	\end{subfigure}
	\hfill
	\begin{subfigure}[b]{0.24\textwidth}
		\includegraphics[scale=0.48]{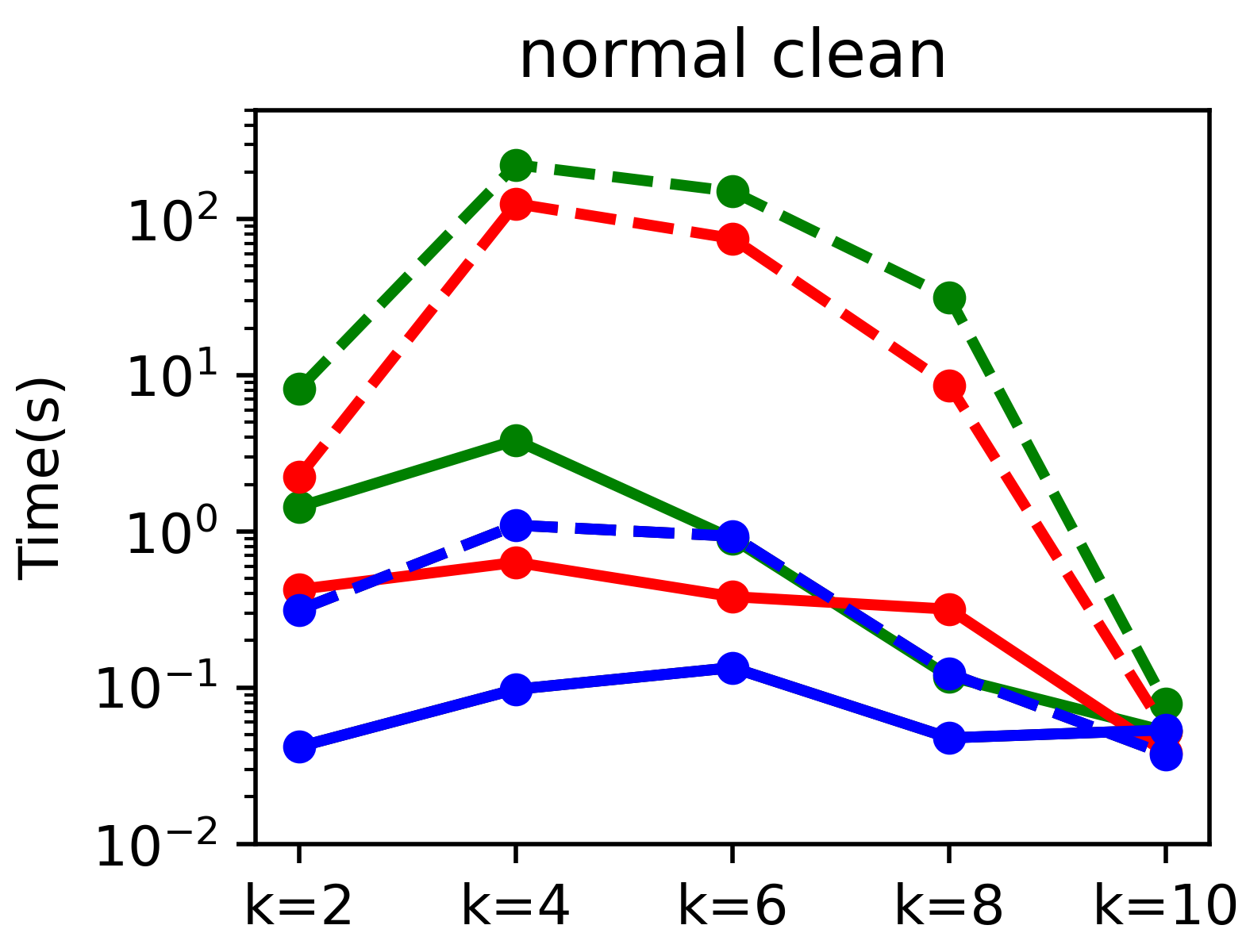}
	\end{subfigure}	
	\hfill	
	\begin{subfigure}[b]{0.24\textwidth}
		\includegraphics[scale=0.48]{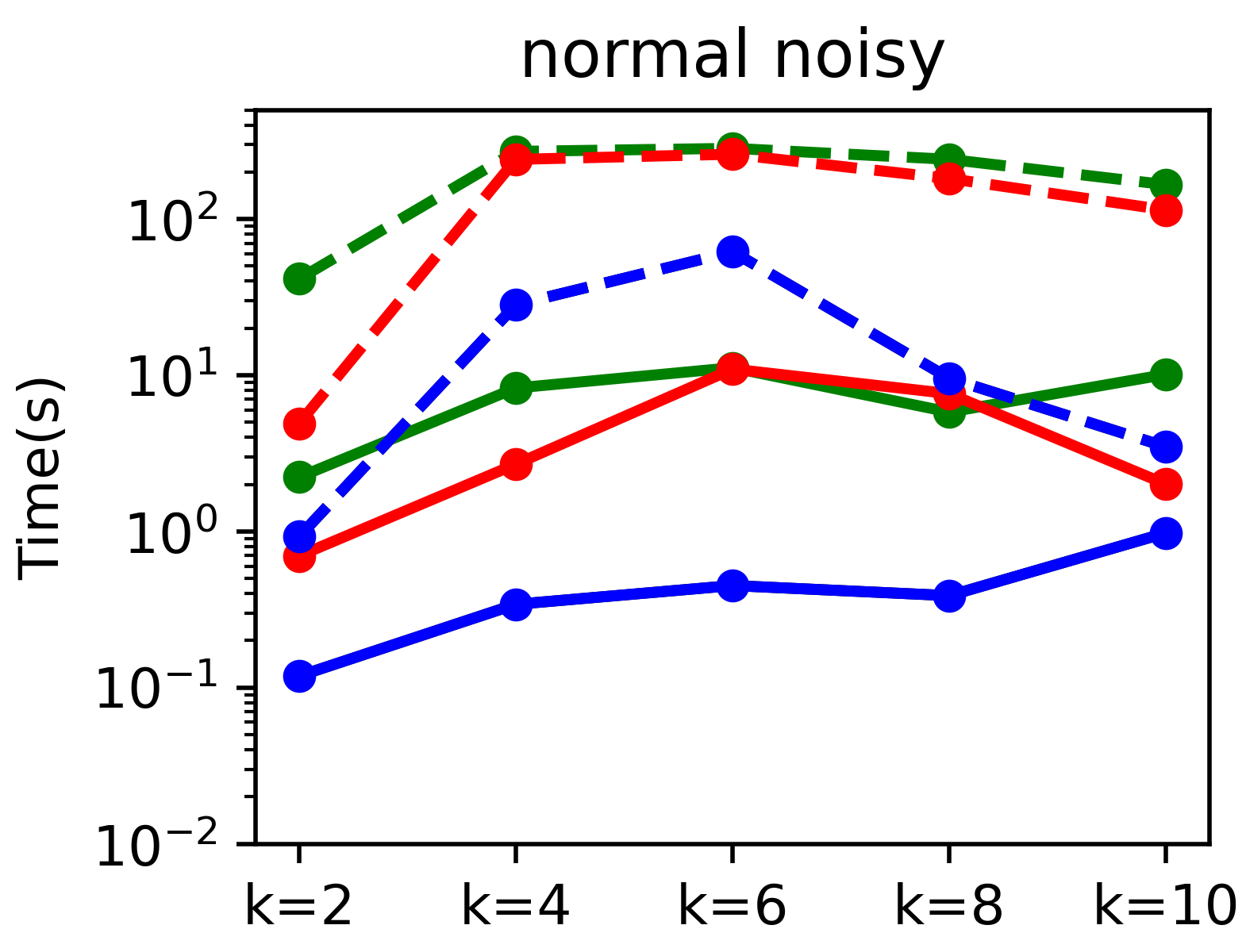}
	\end{subfigure}
\end{figure}

Figure \ref{figure_MLP1_MIPs_time} shows the time taken to solve the master IPs on columns generated by MLP(1).
We observe that MIP(1) takes notably faster to solve than $\text{MIP}_{F}$ and  on most normal-noisy matrices $\text{MIP}_{F}$ runs out of the time budget of $300$ seconds.
Solving both master IPs on columns by MLP($\frac{1}{k}$) also shows us that while solving MIP(1) over a larger set of columns adds only a few seconds for most instances, $\text{MIP}_{F}$  runs out of the time budget of $300$ secs in about half the cases, see Appendix Table \ref{table_MIP1_MIPexact_from_MLP1k}. 
These observations suggest using MIP(1) to find integer feasible solutions in the future as the solution quality is extremely close to that of $\text{MIP}_{F}$ but at a fraction of computational effort. 

\subsection{Accuracy and speed of the IP Formulations.}

In this section we computationally compare the integer programs introduced in Section \ref{section_formulations} and \ref{section_new_objective}. 
CIP due to its polynomial size can be directly given to a general purpose IP solver like CPLEX and we set a time limit of 600 seconds on its running time.
We expect solution times for CIP to grow proportional to $k$ and density of $X$ according to Proposition \ref{lemma_compact_fromulation_lp_0_obj}.
Similarly, we may try to attack the exponential formulation EIP directly by CPLEX. Since however EIP requires the complete enumeration of $2^m$ binary vectors for an input matrix $X$ of size $n\times m$ we can only solve its root LP under 600 seconds in a very few cases. For these few cases however, we observe the objective value of ELP to agree with MLP($\frac{1}{k}$), which gives an experimental confirmation of Proposition \ref{lemma_LPexp_equiv_to_MLP_F}.
In the following experiments, formulation $\text{MIP}_{\text{F}}$ is used on columns generated by MLP($\frac{1}{k}$), while MIP(1) on columns by MLP(1). The final solution of MIP(1) is evaluated under the original $\|\cdot \|_F^2$ objective and that error is reported. As previously, the master LPs are solved with a time limit of $600$ seconds and the master IPs with an additional time limit of $300$ seconds. 

Table \ref{table_MIPexact_MIP1_CIPexact_error} shows the factorisations error in $\|\cdot \|_F^2$ obtained by $\text{MIP}_{\text{F}}$, MIP(1) and CIP and Table \ref{table_MIPexact_MIP1_CIPexact_time} shows the corresponding solution times in seconds.
Each row of Table \ref{table_MIPexact_MIP1_CIPexact_error} and \ref{table_MIPexact_MIP1_CIPexact_time}  corresponds to the average of 10 synthetic matrices of the same size, sparsity and noise. The lowest error results are indicated in boldface. 
We observe that MIP(1) provides the lowest error factorisation in most cases, but CIP gives the lowest error when only looking at $k=2$. The significantly higher error values of $\text{MIP}_{\text{F}}$ are due to the lower quality columns generated by MLP($\frac{1}{k}$) on which it is solved. We emphasise that we do not do branch-and-price when solving MIP(1) or $\text{MIP}_{\text{F}}$. 
Table \ref{table_MIPexact_MIP1_CIPexact_time} shows that MIP(1) is  the fastest in all cases, while CIP runs out of its time limit on all noisy instances for $k=5,10$.
In conclusion, CIP provides very accurate solutions for $k=2$  but it is slower to solve than MIP(1), while for larger $k$'s MIP(1) dominates in both accuracy and speed.
\begin{table}[htbp]
	\centering
	\caption{\label{table_MIPexact_MIP1_CIPexact_error}Factorisation error in $\|\cdot \|_F^2$ of solutions obtained via formulations $\text{MIP}_{\text{F}}$, MIP(1) and $\text{CIP}_{\text{F}}$}
	\begin{tabular}{lccccccccc}
		data  & \multicolumn{3}{c}{k=2} & \multicolumn{3}{c}{k=5} & \multicolumn{3}{c}{k=10} \\
		(n-sparsity-noise) & $\text{MIP}_{\text{F}}$ & MIP(1) & CIP & $\text{MIP}_{\text{F}}$ & MIP(1) & CIP & $\text{MIP}_{\text{F}}$ & MIP(1) & CIP \\
		\hline
\up\down
		20-sparse-clean & 49.6  & \textbf{47.4} & \textbf{47.4} & 20.8  & \textbf{16.6} & 16.7  & \textbf{0.0} & \textbf{0.0} & \textbf{0.0} \\
		20-sparse-noisy & 64.0  & 59.5  & \textbf{59.3} & 42.6  & \textbf{30.3} & 30.7  & 11.2  & \textbf{10.2} & 10.3 \\
		20-normal-clean & 75.0  & 70.0  & \textbf{68.7} & 30.6  & 27.7  & \textbf{26.5} & 0.3   & 0.3   & \textbf{0.0} \\
		20-normal-noisy & 84.6  & 78.9  & \textbf{77.2} & 47.3  & 40.2  & \textbf{40.1} & 11.2  & \textbf{10.7} & 11.2 \\
		\hline
\up\down
		35-sparse-clean & 90.9  & \textbf{84.7} & \textbf{84.7} & 39.1  & \textbf{34.5} & 34.9  & 0.1   & \textbf{0.0} & \textbf{0.0} \\
		35-sparse-noisy & 113.4 & 107.5 & \textbf{106.9} & 84.4  & \textbf{60.5} & 61.7  & 28.4  & \textbf{23.3} & 27.1 \\
		35-normal-clean & 134.2 & 125.0 & \textbf{121.7} & 64.5  & 54.1  & \textbf{53.4} & \textbf{0.0}   & \textbf{0.0} & \textbf{0.0} \\
		35-normal-noisy & 153.6 & 143.1 & \textbf{139.1} & 101.7 & \textbf{80.3} & 81.7  & 31.1  & \textbf{25.5} & 31.1 \\
		\hline
\up\down
		50-sparse-clean & 136.0 & 126.1 & \textbf{125.6} & 61.4  & \textbf{50.6} & 51.5  & 0.1   & \textbf{0.0} & \textbf{0.0} \\
		50-sparse-noisy & 166.2 & \textbf{156.5} & 156.7 & 135.0 & \textbf{89.8} & 93.9  & 49.6  & \textbf{36.7} & 41.4 \\
		50-normal-clean & 215.1 & 198.0 & \textbf{194.3} & 106.1 & \textbf{91.0} & 95.0  & \textbf{0.0}   & \textbf{0.0} & \textbf{0.0} \\
		50-normal-noisy & 237.2 & 218.6 & \textbf{214.2} & 168.6 & 123.9 & \textbf{123.4} & 62.2  & \textbf{44.3} & 61.3 \\
		\hline
	\end{tabular}%
	\bigskip
	\centering
	\caption{  \label{table_MIPexact_MIP1_CIPexact_time}Time in seconds to obtain solutions in Table \ref{table_MIPexact_MIP1_CIPexact_error} via formulations $\text{MIP}_{\text{F}}$, MIP(1) and $\text{CIP}_{\text{F}}$}
	\begin{tabular}{lccccccccc}
		data  & \multicolumn{3}{c}{k=2} & \multicolumn{3}{c}{k=5} & \multicolumn{3}{c}{k=10} \\
		(n-sparsity-noise) & $\text{MIP}_{\text{F}}$ & MIP(1) & CIP & $\text{MIP}_{\text{F}}$ & MIP(1) & CIP & $\text{MIP}_{\text{F}}$ & MIP(1) & CIP \\
		\hline
\up\down
		20-sparse-clean & 1.1   & 0.4   & 1.6   & 4.6   & 0.4   & 169.7 & 0.7   & 0.4   & 1.9 \\
		20-sparse-noisy & 2.7   & 0.6   & 21.8  & 233.7 & 0.8   & 601.6 & 10.9  & 1.8   & 602.9 \\
		20-normal-clean & 15.2  & 3.5   & 56.2  & 303.2 & 5.4   & 600.3 & 3.3   & 1.0   & 15.8 \\
		20-normal-noisy & 31.3  & 5.4   & 295.5 & 336.6 & 17.6  & 600.8 & 65.2  & 8.0   & 602.0 \\
		\hline
\up\down
		35-sparse-clean & 4.0   & 0.8   & 17.3  & 108.4 & 0.9   & 449.8 & 1.9   & 0.5   & 5.3 \\
		35-sparse-noisy & 12.1  & 1.9   & 147.8 & 514.0 & 6.4   & 602.3 & 275.1 & 6.8   & 605.2 \\
		35-normal-clean & 76.0  & 14.2  & 188.6 & 378.5 & 21.8  & 600.8 & 23.2  & 1.6   & 80.6 \\
		35-normal-noisy & 195.3 & 31.8  & 589.7 & 739.3 & 132.1 & 600.7 & 394.7 & 45.3  & 602.4 \\
		\hline
\up\down
		50-sparse-clean & 2.6   & 0.6   & 21.9  & 176.3 & 1.1   & 519.9 & 3.8   & 0.7   & 12.9 \\
		50-sparse-noisy & 28.1  & 2.2   & 285.4 & 827.7 & 6.6   & 602.3 & 523.9 & 6.9   & 605.1 \\
		50-normal-clean & 362.0 & 46.8  & 509.9 & 692.1 & 153.6 & 602.1 & 187.2 & 2.5   & 139.4 \\
		50-normal-noisy & 601.6 & 194.8 & 578.2 & 903.9 & 341.1 & 601.0 & 649.8 & 146.2 & 601.6 \\
		\hline
\up\down
	\end{tabular}%
\end{table}%

\subsection{Binary matrix completion.}

In this section we explore how successful our approach is at recovering missing entries of incomplete binary matrices. 
We create an incomplete dataset of our synthetic matrices by deleting $5,10,\dots, 30\%$ of the entries of each matrix.
This way, after computing a rank-$k$ factorisation of the incomplete matrix, we can easily compare to the corresponding original matrix to see how many of the entries we have recovered successfully. 
Since our synthetic matrices are generated to be of Boolean rank at most $10$, we cannot expect to recover all the entries by a rank-$k$ completion with $k<10$ and thus we perform the experiments with $k=10$.

\begin{figure}[htbp]
	\caption{ {Rank-$10$ binary matrix completion of artificial matrices with $5-30\%$ missing entries }}
	\centering
	\includegraphics[scale=0.5]{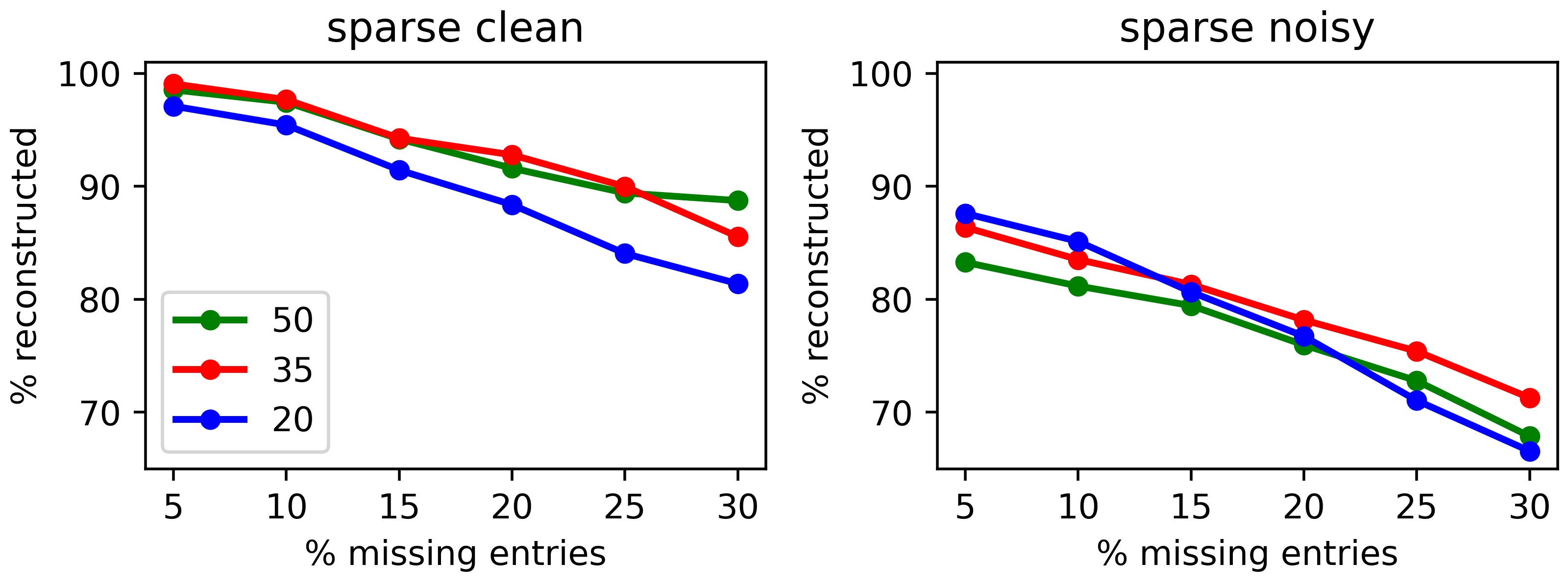}
	\includegraphics[scale=0.5]{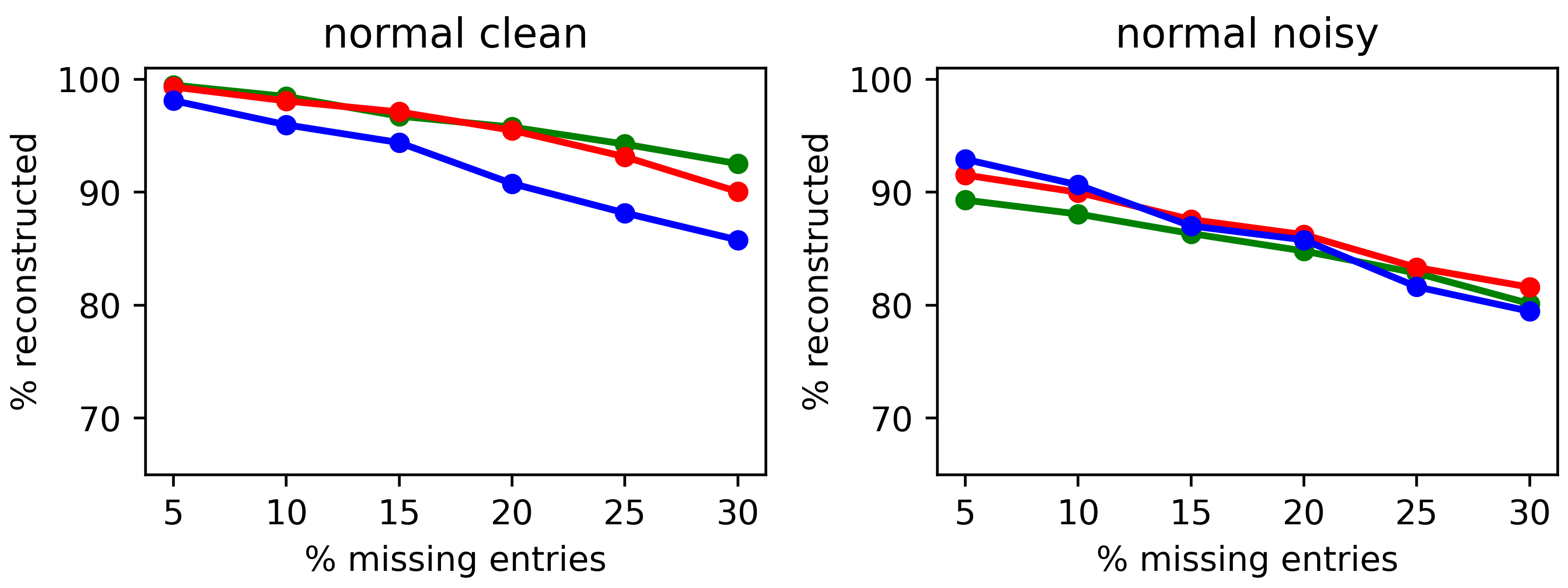}
	\label{figure_completion}
\end{figure}

Figure \ref{figure_completion} shows the reconstruction percentage against the percentage of missing entries when solving MIP(1) on columns generated by MLP(1) on the incomplete matrices. As previously, the three colours correspond to dimensions of the matrices: green to $50 \times 20$, red to $35 \times 20$ and blue to $20 \times 20$.
We define the percentage of reconstruction  as $100*(1-\|X-A\circ B\|_F^2/ \|X \|_F^2)$ where $X$ is the original complete matrix and $A\circ B$ is the rank-$k$ factorisation of the incomplete matrix.
As expected the recovery percentage decreases with the percentage of missing entries and 
clean matrices are better recovered than noisy ones. All in all, we see a very high percentage of the entries can be recovered by MIP(1).

\subsection{Comparing integer programming approaches against heuristics.} \label{sec-comp-heur}
In this section, we compare our integer programming approaches against the most widely used $k$-BMF heuristics on real-world datasets.
The heuristic algorithms we evaluate include the ASSO algorithm \cite{Miettinen:2006:PKDD, Miettinen:2008:DBP:1442800.1442809}, 
the alternating iterative local search algorithm (ASSO++) of \cite{Barahona:2019}  which uses ASSO as a starting point, 
and the penalty objective formulation (pymf) of \cite{Zhang:2007} via the implementation of \cite{Schinnerl:2017}. We also compute rank-$k$ NMF and binarise it with a threshold of $0.5$. 
The exact details and parameters used in the computations can be found in Appendix \ref{appendix_k_BMF_heuristics}. 
In addition, we use a new heuristic  which sequentially finds $k$ rank-$1$ binary matrices using any heuristic for Bipartite Binary Quadratic Programming as a subroutine. We refer to this heuristic outlined in Algorithm \ref{greedy_rank_k} as \textit{$k$-Greedy} as the subroutine we use to compute the rank-$1$ binary matrices is the greedy algorithm of \cite{Punnen:2012}. 
\begin{algorithm}[htbp]
\SetAlgoNoLine
	\caption{\label{greedy_rank_k}Greedy algorithm for $k$-BMF ($k$-Greedy)}
	Input: $X\in\B^{n\times m}$, $k\in \Z_+$. \\
	Set $H \in \{-1,0,1 \}^{n\times m}$ to $h_{ij}=2 x_{ij}-1$ for $(i,j)\in E\cup \overline{E}$ and $h_{ij}=0$ otherwise.\\
	\For{$\ell\in[k]$}{
		$\mvec{a}, \mvec{b} = \text{BBQP} (H)$ \tcp{compute a rank-1 binary matrix via any algorithm for BBQP}
		$A_{:,\ell} = \mvec{a}$ \\
		$B_{\ell,:} = \mvec{b}^\top$\\
		$H[\mvec{a} \mvec{b}^\top==1] = 0$ \tcp{set entries of $H$ to zero that are covered}
	}
	Output: $A\in\B^{n\times k}$, $B\in\B^{k\times m}$
\end{algorithm}

We solve CIP using CPLEX with a time limit of $20$ mins and provide the heuristic solution of $k$-Greedy as a warm start to it.
The column generation approach results are obtained by generating columns for $20$ mins using formulation MLP(1) with a warm start of initial rank-1 binary matrices obtained from $k$-Greedy, then solving MIP(1) over the generated columns with a time limit of $10$ mins.
Table \ref{COMP} shows the factorisation error in $\|\cdot\|_F^2$ after evaluating the above described methods on all real-world datasets without missing entries for $k=2,5,10$. 
The best result for each instance is indicated in boldface. 
We observe that CG provides the strictly smallest error for 8 out of 12 instances.
\begin{table}[htbp]
	\centering
	\caption{Comparison of factorisation error in $\|\cdot\|_F^2$ for two IP based methods and five $k$-BMF heuristics}
	\begin{tabular}{clrrrrrrr}
		&       & \multicolumn{1}{l}{MIP(1)} & \multicolumn{1}{l}{CIP} & \multicolumn{1}{l}{ASSO++} & \multicolumn{1}{l}{k-Greedy} & \multicolumn{1}{l}{pymf} & \multicolumn{1}{l}{ASSO} & \multicolumn{1}{l}{NMF} \\
		\hline
\up\down
		\multirow{4}[2]{*}{k=2} & zoo   & 272   & \textbf{271} & 276   & 323   & 274   & 367   & 295 \\
		& heart & \textbf{1185} & 1187  & 1187  & 1187  & 1241  & 1251  & 1273 \\
		& lymp  & 1192  & \textbf{1184} & 1202  & 1201  & 1225  & 1352  & 1427 \\
		& apb   & \textbf{776} & \textbf{776} & \textbf{776} & \textbf{776} & 794   & 778   & 820 \\
		\hline
\up\down
		\multirow{4}[2]{*}{k=5} & zoo   & \textbf{126} & 129   & 133   & 218   & 153   & 354   & 135 \\
		& heart & \textbf{737} & 738   & 738   & 738   & 813   & 887   & 1190 \\
		& lymp  & \textbf{982} & 1026  & 1039  & 1053  & 1067  & 1484  & 1112 \\
		& apb   & \textbf{684} & 688   & 694   & 688   & 733   & 719   & 729 \\
		\hline
\up\down
		\multirow{4}[2]{*}{k=10} & zoo   & \textbf{39} & 72    & 55    & 175   & 80    & 377   & 319 \\
		& heart & 425   & 529   & \textbf{419} & 565   & 483   & 694   & 896 \\
		& lymp  & \textbf{728} & 829   & 812   & 859   & 952   & 1525  & 1102 \\
		& apb   & \textbf{573} & 605   & 591   & 606   & 611   & 661   & 660 \\
		\hline
	\end{tabular}%
	\label{COMP}
\end{table}%

While integer programming based approaches are able to handle missing entries by simply setting the objective coefficients of the missing entries to $0$, the $k$-BMF heuristics ASSO, ASSO++ and pymf cannot so simply be adjusted. 
Non-negative matrix factorisation however, has an available implementation that can handle missing entries \cite{Li:2012, Li:Code:2012}. Our next experiment compares our integer programming approaches against $k$-Greedy and NMF on the real datasets that have missing entries. Table \ref{table_real_missing} shows the results with the lowest error results indicated in boldface. For $k=2$, $k$-Greedy provides very accurate solutions which MIP(1) and CIP fail to improve on in $3$ out of $4$ instances. For $k=5,10$ however, MIP(1) produces notably lower error factorisations than the other methods.
\begin{table}[htbp]
	\centering
	\caption{Comparison of factorisation error in $\|\cdot\|_F^2$ for real-world data with missing entries}
	\begin{tabular}{clrrrr}
		&       & \multicolumn{1}{l}{MIP(1)} & \multicolumn{1}{l}{CIP} & \multicolumn{1}{l}{k-Greedy} & \multicolumn{1}{l}{NMF} \\
		\hline
\up\down
		\multirow{4}[2]{*}{k=2} & tumor & \textbf{1352} & \textbf{1352} & \textbf{1352} & 1792 \\
		& hepatitis & \textbf{1264} & 1344  & 1416  & 1346 \\
		& audio & \textbf{1419} & \textbf{1419} & \textbf{1419} & 2361 \\
		& votes & \textbf{1246} & \textbf{1246} & \textbf{1246} & 1268 \\
		\hline
		\up\down
		\multirow{4}[2]{*}{k=5} & tumor & \textbf{962} & 993   & 1004  & 1832 \\
		& hepatitis & \textbf{1138} & 1229  & 1238  & 1618 \\
		& audio & \textbf{1064} & 1078  & 1094  & 2361 \\
		& votes & \textbf{779} & 853   & 853   & 2353 \\
		\hline
		\up\down
		\multirow{4}[2]{*}{k=10} & tumor & \textbf{514} & 632   & 646   & 1949 \\
		& hepatitis & \textbf{907} & 1048  & 1056  & 2159 \\
		& audio & \textbf{765} & 881   & 881   & 2361 \\
		& votes & \textbf{240} & 701   & 706   & 3189 \\
		\hline
	\end{tabular}%
	\label{table_real_missing}%
\end{table}%
\section{Conclusions and further work.}
In this paper we investigated the rank-$k$ binary matrix factorisation problem from an integer programming perspective. We analysed a compact and two exponential size integer programming formulations for the problem and made a comparison on the strength of the formulations' LP-relaxations. We introduced a new objective function, which slightly differs from the traditional squared Frobenius objective in attributing a weight to zero entries of the input matrix that is proportional to the number of times the zero is erroneously covered in a rank-$k$ factorisation. In addition, we discussed a computational approach based on column generation to solve one of the exponential size formulations and reported several computational experiments to demonstrate the applicability of our formulations on real world and artificial datasets. 

Future research directions that could be explored include developing faster exact algorithms for the pricing problem and once the pricing problems are solved more efficiently, a full branch-and-price implementation would be interesting to explore.

%
%
%
 \begin{APPENDICES}
 \section{Heuristics for the pricing problem.}
 \label{appendix_pricing}

 The greedy algorithm of \cite{Punnen:2012} to solve the Bipartite Binary Quadratic Program in Equation \eqref{equation_PP_BBQ} aims to set entries of $\mvec{a}$ and $\mvec{b}$ to $1$ which correspond to rows and columns of $H$ with the largest positive weights. 
 In the first phase  of the algorithm, the row indices $i$ of $H$ are put in decreasing order according to their sum of positive entries, so $\gamma^+_{i} \ge \gamma^+_{i+1}$ where $ \gamma^+_i:=\sum_{j=1}^m \max(0,h_{ij})$.  
 Then sequentially according to this ordering, $a_i$ is set to $1$ if $\sum_{j=1}^m \max ( 0, \sum_{\ell=1}^{i-1}  a_\ell h_{\ell j} ) < \sum_{j=1}^m \max ( 0, \sum_{\ell=1}^{i} a_\ell h_{\ell j})$ and $0$ otherwise. In the second phase, $b_j$ is set to $1$ if $(\mvec{a}^\top H)_j >0$, $0$ otherwise. An efficient implementation of the greedy algorithm due to \cite{Punnen:2012} is given in Algorithm \ref{alg1}. 
 
 \begin{algorithm}[htbp]
 	\caption{Greedy Algorithm for BBQP}
 	\label{alg1}
 	\begin{multicols}{2}
\SetAlgoNoLine
 		Input: $H\in\R^{n\times m}$ \\
 		\textbf{Phase I.} Order $i\in[n]$ so that $\gamma^+_{i} \ge \gamma^+_{i+1}$.\\
 		Set $\mvec{a}=\mvec{0}_n$, $\mvec{s}=\mvec{0}_m$. \\
 		\For{$i\in[n]$}{
 			$f_{0} = \sum_{j=1}^m \max ( 0, s_j)$\\
 			$ f_1= \sum_{j=1}^m \max ( 0, s_j+ h_{ij})$\\
 			\If{$f_{0} < f_1$}{ 
 				Set $a_i = 1$, $\mvec{s} = \mvec{s} + \mvec{h}_i$
 			}
 		}
 		\columnbreak
 		$\quad$\\
 		\textbf{Phase II.} \\
 		Set $\mvec{b}=\mvec{0}_m$.\\
 		\For{$j\in[m]$}{
 			\If{$(\mvec{a}^\top H)_j>0$}{
 				Set $b_j = 1$}
 		}
 		Output: $\mvec{a}\in \B^n, \mvec{b}\in\B^m$
 	\end{multicols}
 \end{algorithm}
 
 There are many variants of Algorithm \ref{alg1} one can explore.
 First, the solution greatly depends on the ordering of $i$'s in the first phase. If for some $i_1\not = i_2$ we have $\gamma^{+}_{i_1} = \gamma^{+}_{i_2}$, comparing the sum of negative entries of rows $i_1$ and $i_2$ can put more ``influential'' rows of $H$ ahead in the ordering. Let us call this ordering the \textit{revised ordering} and the one which only compares the positive sums as the \textit{original ordering}.
 Another option is to use a completely \textit{random order} of $i$'s or to apply a small perturbation to sums $\gamma^+_i$ to get a \textit{perturbed} version of the revised or original ordering. None of the above ordering strategies clearly dominates the others in all cases but they are fast to compute hence one can evaluate all five ordering strategies (original, revised, original perturbed, revised perturbed, random) and pick the best one.
 Second, the algorithm as presented above first fixes $\mvec{a}$ and then $\mvec{b}$. Changing the order of fixing $\mvec{a}$ and $\mvec{b}$ can yield a different result hence it is best to try for both $H$ and $H^\top$. In general, it is recommended to start the first phase on the smaller dimension \cite{Punnen:2012}. 
 Third, the solution from Algorithm \ref{alg1} may be improved by computing the optimal $\mvec{a}$ with respect to fixed $\mvec{b}$. This idea then can be used to fix $\mvec{a}$ and $\mvec{b}$ in an alternating fashion and stop when no changes occur in either. We summarise this alternating heuristic in Algorithm \ref{alg_iterating}
 
 \begin{algorithm}[htbp]
 	\caption{Alternating Heuristic for BBQP}
 	\label{alg_iterating}
 	Input: $H\in\R^{n\times m}$, $\mvec{a}^{(0)}\in\B^n,\mvec{b}^{(0)}\in\B^m$. \\
 	\For{$\ell=1,2,\dots$}{
 		$\mvec{a}^{(\ell)} [ H\mvec{b}^{(\ell-1)} >0] = 1$\\
 		$\mvec{a}^{(\ell)} [ H\mvec{b}^{(\ell-1)} \le 0] = 0$\\		
 		\If{$\mvec{a}^{(\ell)} == \mvec{a}^{(\ell-1)} $}{Break}
 		$\mvec{b}^{(\ell)} [ (\mvec{a}^{(\ell)})^\top H >0] = 1$\\
 		$\mvec{b}^{(\ell)} [ (\mvec{a}^{(\ell)})^\top H \le0] = 0$\\
 		\If{$\mvec{b}^{(\ell)} == \mvec{b}^{(\ell-1)} $}{Break}
 	}
 	Output: $\mvec{a}^{(\ell)}\in\B^n,\mvec{b}^{(\ell)}\in\B^m$
 \end{algorithm}
 In Section \ref{subsection_cg_test} we use the above described heuristics for the pricing problem in column generation. At each iteration of the column generation procedure, 30 variants of Algorithm \ref{alg1} are computed to obtain an initial feasible solution to PP. The 30 variants of the greedy algorithm use the original and revised ordering, their transpose and perturbed version and 22 random orderings. All greedy solutions are improved by the alternating heuristic until no further improvement is found.

 \section{Synthetic data.}
 \label{appendix_synthetic_data}
 Table \ref{table_synthetic_data} gives a summary of the parameters used to generate our synthetic dataset.
 For a synthetic binary matrix $X$, $n\times m$ is the dimension of $X$, $\kappa$ is the Boolean rank which was used to generate $X$, and $n'\times m'$ is the dimension obtained after removing zero and duplicate row and columns of $X$.
 \begin{table}[htbp]
 	\centering
 	\caption{Parameters of the synthetic dataset}
 	\begin{tabular}{lcccc c c}
 		(n-sparsity-noise) & $n \times m$  & \multicolumn{1}{l}{$\kappa$} & \multicolumn{1}{l}{0s\%} & \multicolumn{1}{l}{noise\%}  & \#instances & $n'\times m'$ \\
		\hline
\up\down
 		20-sparse-clean & \multirow{4}[2]{*}{20 $\times$ 20} & \multirow{4}[2]{*}{10} & \multirow{2}[1]{*}{75} & 0 & \multirow{4}[2]{*}{10} & $14 \times 15$\\
 		20-sparse-noisy &       &       &       & 5 &  & $19 \times 19$\\
 		20-normal-clean &       &       & \multirow{2}[1]{*}{50} & 0 && $18 \times 18$\\
 		20-normal-noisy &       &       &       & 5 & & $19 \times 20$\\
		\hline
\up\down
 		35-sparse-clean & \multirow{4}[2]{*}{35 $\times$  20} & \multirow{4}[2]{*}{10} & \multirow{2}[1]{*}{75} & 0 & \multirow{4}[2]{*}{10} & $22 \times 15$\\
 		35-sparse-noisy &       &       &       & 5 &   & $31 \times 19$\\
 		35-normal-clean &       &       & \multirow{2}[1]{*}{50} & 0  & & $29 \times 18$\\
 		35-normal-noisy &       &       &       & 5 & & $34 \times 20$\\
		\hline
\up\down
 		50-sparse-clean & \multirow{4}[2]{*}{50 $\times$  20} & \multirow{4}[2]{*}{10} & \multirow{2}[1]{*}{75} & 0 &  \multirow{4}[2]{*}{10} & $30 \times 15$\\
 		50-sparse-noisy &       &       &       & 5 && $45 \times 20$\\
 		50-normal-clean &       &       & \multirow{2}[1]{*}{50} & 0 && $40 \times 18$\\
 		50-normal-noisy &       &       &       & 5& & $48 \times 20$\\
		\hline
 	\end{tabular}%
 	\label{table_synthetic_data}%
 \end{table}%
 
 \section{Real world data.}
 \label{appendix_real_data}
 
 The following datasets were used in the experiments:
 \begin{itemize}
 	\item 
 	The Zoo dataset (\textit{zoo})  \cite{ZOO} describes $101$ animals with $16$ characteristic features. All but one feature is binary. The categorical column which records the number of legs an animal has, is converted into two new binary columns indicating if the number of legs is \textit{less than or equal} or \textit{greater} than four. The size of the resulting fully binary matrix is $101 \times 17$.
 	
 	\item 
 	The Primary Tumor dataset (\textit{tumor}) \cite{TUMOR} contains observations on $17$ tumour features detected in $339$ patients. The features are represented by $13$ binary variables and $4$ categorical variables with discrete options. The $4$ categorical variables are converted into $11$ binary variables representing each discrete option. Two missing values in the binary columns are left as missing values. The final dimension of the binary matrix  is $339 \times 24$ with 670 missing values.
 	
 	\item 
 	The Hepatitis dataset (\textit{hepat}) \cite{HEP}  consists of 155 samples of medical data of patients with hepatitis. The 19 features of the dataset can be used to predict whether a patient with hepatitis will live or die. 
 	6 of the 19 features take numerical values and are converted into 12 binary features corresponding to options: \textit{less than or equal to the median value}, and \textit{greater than the median value}. The column that stores the sex of patients is converted into two binary columns corresponding to labels man and female. The remaining 12 columns take values \textit{yes} and \textit{no} and are converted into 24 binary columns. The missing values in the raw dataset are left as missing in the binary dataset as well.
 	The final dimension of the binary matrix is $155 \times 38$ with $334$ missing values.
 	
 	\item 
 	The SPECT Heart dataset (\textit{heart}) \cite{SPECT} describes cardiac Single Proton Emission Computed Tomography images of $267$ patients by $22$ binary feature patterns. $25$ patients' images contain none of the features and are dropped from the dataset, hence the final dimension of the binary matrix is $242 \times 22$.
 	
 	\item The Lymphography dataset (\textit{lymp})  \cite{LIM} contains data about lymphography examination of $148$ patients. $8$ features take categorical values and are expanded into $33$ binary features representing each categorical value. One column is numerical and we convert it into two binary columns corresponding to options: \textit{less than or equal to median value}, and \textit{larger than median value}. The final dimension of the fully binary matrix is $148\times 44$.
 	
 	\item
 	The Audiology Standardized dataset (\textit{audio}) \cite{AUDIO} contains clinical audiology records on $226$ patients. The $69$ features include patient-reported symptoms, patient history information, and the results of routine tests which are needed for the evaluation and diagnosis of hearing disorders. $9$ features that are categorical valued are binarised into $34$ new binary variables indicating if a discrete option is selected.  The missing values in the raw dataset are left as missing in the binary dataset as well. The final dimension of the binary matrix is $226\times 92$ with $899$ missing values.
 	
 	\item 
 	The Amazon Political Books dataset (\textit{books}) \cite{BOOKS} contains binary data about $105$ US politics books sold by Amazon.com. Columns correspond to books and rows represent frequent co-purchasing of books by the same buyers. The dimension of the binary matrix is $105\times 105$.
 	
 	\item 
 	The 1984 United States Congressional Voting Records dataset (\textit{votes})\cite{VOTING} includes votes for each of the U.S. House of Representatives Congressmen on the $16$ key votes identified by the CQA. The $16$ categorical variables taking values of ``voted for'', ``voted against'' or ``did not vote'', are converted into $16$ binary features taking value $1$ for ``voted for'', value $0$ for ``voted against'' and a missing value indicates ``did not vote''.  The final dimension of the binary matrix is $435\times 16$ with $392$ missing values.
 	
 \end{itemize}

 \section{Obtaining integer feasible solutions.}
 \label{appendix_obtaining_integer_sols}
 In this section we give additional numerical results supporting our conclusions drawn in Section \ref{subsection_obtaining_integer_sols}.
 Table \ref{table_MIP1_MIPexact_from_MLP1} shows the factorisation error measured in $\|\cdot\|_F^2$ of integer feasible solutions obtained by solving MIP(1) and $\text{MIP}_{F}$ over columns generated by MLP(1).
 MIP(1) takes significantly faster to solve than $\text{MIP}_{F}$ but the absolute difference in error between solutions produced by MIP(1) and $\text{MIP}_{F}$ is at most $1$, except for the last row in column $k=5$ where $\text{MIP}_{F}$ runs out of the time budget of $300$ seconds and produces higher error solutions than MIP(1).
 
 Table \ref{table_MIP1_MIPexact_from_MLP1k} shows the result of an analogous experiment where the columns used are generated by MLP($\frac{1}{k}$). 
 Since MLP($\frac{1}{k}$) is slower to solve than MLP(1), more columns are generated during CG and the master IPs have a harder task on selecting $k$ columns from a larger set of columns in Table \ref{table_MIP1_MIPexact_from_MLP1k}. However, while solving MIP(1) over a larger set of columns adds only a few seconds for most instances, $\text{MIP}_{F}$  runs out of the time budget of $300$ secs in about half the cases. This is also demonstrated  in the error difference, with solutions by MIP(1) having smaller error than solutions by $\text{MIP}_{F}$ in most cases. 
 
 \begin{table}[htbp]
 	\centering
 	\caption{Error in $\|\cdot\|_F^2$ (and runtime in seconds) of integer solutions by MIP(1) and $\text{MIP}_{\text{F}}$ on columns by MLP(1)}
 	\begin{tabular}{lcccccc}
 		data  & \multicolumn{2}{c}{k=2} & \multicolumn{2}{c}{k=5} & \multicolumn{2}{c}{k=10} \\
 		(n-sparsity-noise) & MIP(1) & $\text{MIP}_{\text{F}}$ & MIP(1) & $\text{MIP}_{\text{F}}$ & MIP(1) & $\text{MIP}_{\text{F}}$ \\
		\hline
\up\down
 		20-sparse-clean & 47 (0.0) & 47 (0.0) & 16 (0.0) & 16 (0.0) & 0 (0.0) & 0 (0.0) \\
 		20-sparse-noisy & 59 (0.0) & 59 (0.0) & 30 (0.0) & 30 (0.0) & 10 (0.0) & 10 (0.0) \\
 		20-normal-clean & {70 (0.0)} & \textbf{69} (0.3) & 27 (0.1) & 27 (2.7) & 0 (0.0) & 0 (0.0) \\
 		20-normal-noisy & 78 (0.1) & 78 (0.9) & {40 (0.5)} & \textbf{39} (76.5) & 10 (0.5) & 10 (3.4) \\
		\hline
\up\down
 		35-sparse-clean & 84 (0.0) & 84 (0.1) & 34 (0.0) & 34 (0.1) & 0 (0.0) & 0 (0.0) \\
 		35-sparse-noisy & 107 (0.0) & 107 (0.1) & 60 (0.0) & 60 (0.6) & 23 (0.1) & 23 (0.2) \\
 		35-normal-clean & {125 (0.4)} & \textbf{124} (2.2) & {54 (0.8)} & \textbf{53} (154.8) & 0 (0.0) & 0 (0.1) \\
 		35-normal-noisy & {143 (0.6)} & \textbf{141} (4.9) & 80 (4.1) & 80 (245.4) & {25 (2.0)} & \textbf{24} (114.2) \\
		\hline
\up\down
 		50-sparse-clean & 126 (0.0) & 126 (0.0) & 50 (0.0) & 50 (0.1) & 0 (0.0) & 0 (0.0) \\
 		50-sparse-noisy & 156 (0.0) & 156 (0.1) & 89 (0.0) & 89 (0.2) & 36 (0.0) & 36 (0.2) \\
 		50-normal-clean & {198 (1.4)} & \textbf{197} (8.2) & 91 (30.9) & 91 (173.4) & 0 (0.1) & 0 (0.1) \\
 		50-normal-noisy & 218 (2.2) & 218 (41.4) & \textbf{123} (39.7) & {126 (271.1)} & 44 (10.1) & 44 (165.8) \\
		\hline
 	\end{tabular}%
 	\label{table_MIP1_MIPexact_from_MLP1}%
 	\bigskip
 	\centering
 	\caption{Error in $\|\cdot\|_F^2$ (and runtime in seconds) of integer  solutions by MIP(1) and $\text{MIP}_{\text{F}}$ on columns  by MLP($\frac{1}{k}$)}
 	\begin{tabular}{lcccccc}
 		data  & \multicolumn{2}{c}{k=2} & \multicolumn{2}{c}{k=5} & \multicolumn{2}{c}{k=10} \\
 		(n-sparsity-noise) & MIP(1) & $\text{MIP}_{\text{F}}$ & MIP(1) & $\text{MIP}_{\text{F}}$ & MIP(1) & $\text{MIP}_{\text{F}}$ \\
 				\hline
 		\up\down
 		20-sparse-clean & 50 (0.0) & 50 (0.2) & 21 (0.0) & 21 (2.6) & 0 (0.0) & 0 (0.0) \\
 		20-sparse-noisy & 64 (0.0) & 64 (0.6) & \textbf{42} (0.1) &  {43 (219.0)} & 11 (0.2) & 11 (6.3) \\
 		20-normal-clean & {76 (0.2)} & \textbf{75} (3.9) & \textbf{30} (0.5) & {31 (289.6)} & 0 (0.1) & 0 (0.2) \\
 		20-normal-noisy & 85 (0.3) & 85 (6.3) & 47 (1.2) & 47 (300.4) & 11 (0.6) & 11 (54.2) \\
		\hline
\up\down
 		35-sparse-clean & 91 (0.0) & 91 (1.5) & 39 (0.2) & 39 (98.9) & 0 (0.1) & 0 (0.1) \\
 		35-sparse-noisy & 114 (0.1) & 113 (3.1) & \textbf{81} (0.5) & {84 (300.7)} & 28 (0.3) & 28 (229.9) \\
 		35-normal-clean & {136 (1.0)} & \textbf{134} (19.1)& \textbf{61} (2.0) & {65 (300.8)} & 0 (0.8) & 0 (11.9) \\
 		35-normal-noisy & 154 (1.6) & 154 (58.9) & \textbf{93} (6.2) & {102 (301.3)} & \textbf{28} (2.1) & {31 (301.0)} \\
		\hline
\up\down
 		50-sparse-clean & {137 (0.0)} & \textbf{136} (0.8) & 61 (0.2) & 61 (160.0) & 0 (0.8) & 0 (0.2) \\
 		50-sparse-noisy & {167 (0.1)} & \textbf{166} (6.5) & \textbf{128} (0.7) & {135 (301.5)} & \textbf{46} (0.6) & {50 (301.5)} \\
 		50-normal-clean & 215 (2.2) & 215 (131.6) & \textbf{100} (34.4) & {106 (302.1)} & 0 (0.8) & 0 (153.7) \\
 		50-normal-noisy & 238 (5.7) & \textbf{237} (226.4) & \textbf{149} (95.8) & {169 (302.9)} & \textbf{51} (39.4) & {62 (302.5)} \\
		\hline
 	\end{tabular}%
 	\label{table_MIP1_MIPexact_from_MLP1k}%
 \end{table}%

 \section{Heuristics for \texorpdfstring{$k$}{k}-BMF.}
 \label{appendix_k_BMF_heuristics}
 The following methods were evaluated for the comparison in Tables \ref{COMP} and \ref{table_real_missing}.
 \begin{itemize}
 	\item For the alternating iterative local search algorithm  of \cite{Barahona:2019}  (ASSO++) we obtained the code from the author's github page, see the reference. The code implements two variants of the algorithm and we report the smaller error solution from two variants of it. 
 	\item For the method of  \cite{Zhang:2007}, we used a python implementation in the package \texttt{pymf}, see \cite{Schinnerl:2017} and we ran it for 10000 iterations.
 	\item We evaluated the heuristic method ASSO \cite{Miettinen:2006:PKDD} which  depends on a parameter and we report the best results across nine parameter settings ($\tau \in \{0.1, 0.2,\dots ,0.9 \}$). The code was obtained form the webpage of the author: \texttt{https://people.mpi-inf.mpg.de/~pmiettin/src/DBP-progs/}. We observe that ASSO does not return monotone solutions and sometimes we get a higher error solution for a higher value of $k$. 
 	\item  We computed rank-$k$ non-negative matrix factorisation (NMF) and binarise it by a threshold of $0.5$: after an NMF is obtained, values greater than $0.5$ are set to $1$, otherwise to $0$.  For the computation of NMF we used the function \texttt{non\_negative\_factorization} from the \texttt{sklearn.decomposition} module in python when the binary matrix has no missing entries, and for incomplete binary matrices we used the Matlab implementation in \cite{Li:Code:2012, Li:2012}.
 	\item The heuristic $k$-greedy algorithm was ran with $70$ random seeds and the subroutine for BBQP used the greedy and alternating algorithms for BBQP given in Algorithms \ref{alg1}, \ref{alg_iterating}. In addition, the $k$-greedy algorithm can be run on a preprocessed or original matrix and we tried both ways. For each instance the lowest error factorisation is reported.
 \end{itemize}

 \end{APPENDICES}

\section*{Acknowledgments.}
During the completion of this work R.A.K was supported by a doctoral scholarship from The Alan Turing Institute and the Office for National Statistics.


\bibliographystyle{informs2014} 
\bibliography{references.bib} 


\end{document}